\documentclass[twocolumn]{article}
\usepackage[utf8]{inputenc}  
\usepackage[T1]{fontenc}     
\usepackage{textcomp}        
\usepackage{graphicx}
\usepackage{cite}
\usepackage{amsmath, amssymb, amsthm, amsfonts, cases}
\usepackage{mathrsfs}
\usepackage{url}
\usepackage{todo}
\usepackage{authblk}
\usepackage[usenames]{color}
\usepackage{geometry}
\usepackage{indentfirst}
\allowdisplaybreaks[4]

\theoremstyle{plain}
\newtheorem{thm}{Theorem}[section]

\newtheorem{pro}[thm]{Problem}

\newtheorem{prop}[thm]{Proposition}

\theoremstyle{definition}
\newtheorem{defn}{Definition}[section]
\newtheorem{ass}{Assumption}[section]
\newtheorem{rmk}{Remark}[section]

\makeatletter\@addtoreset{equation}{section} \makeatother

\begin{document}		
	\title{ Discrete-Time LQ Stochastic Two Person Nonzero Sum Difference Games With Random Coefficients:~Closed-Loop Nash Equilibrium
		\thanks{Qingxin Meng was supported by the Key Projects of Natural Science Foundation of Zhejiang Province (No. LZ22A010005) and the National
			Natural Science Foundation of China ( No.12271158).}}

	\date{}
	
	\author[a]{Qingxin Meng}
	\author[b]{YiWei Wu
\footnote{Corresponding author. E-mail: elorywoo@gmail.com (Y. Wu), mqx@zjhu.edu.cn (Q. Meng).}}

	\affil[a]{\small{Department of Mathematical Sciences, Huzhou University, Zhejiang 313000,PR  China}}
	\affil[b]{\small{School of Mathematical Sciences, South China Normal University, Guangzhou 510631, China}}

	\maketitle

	\begin{abstract}
This paper investigates closed-loop Nash equilibria for discrete-time linear-quadratic (LQ) stochastic nonzero-sum difference games with random coefficients. Unlike existing works, we consider randomness in both state dynamics and cost functionals, leading to a complex structure of fully coupled cross-coupled stochastic Riccati equations (CCREs). The key contributions lie in characterizing the equilibrium via state-feedback strategies derived by decoupling stochastic Hamiltonian systems governed by two symmetric CCREs-these random coefficients induce a higher-order nonlinear backward stochastic difference equation (BS$\triangle$E) system, fundamentally differing from deterministic counterparts. Under minimal regularity conditions, we establish necessary and sufficient conditions for closed-loop Nash equilibrium existence, contingent on the regular solvability of CCREs without requiring strong assumptions. Solutions are constructed using a dynamic programming principle (DPP), linking equilibrium strategies to coupled Lyapunov-type equations. Our analysis resolves critical challenges in modeling inherent randomness and provides a unified framework for dynamic decision-making under uncertainty.
	\end{abstract}

	
	\section{ Introduction}
	
	Nash equilibrium and LQ differential games been a significant area of research across various fields, including engineering, economics, and biology. Nash equilibrium, initially introduced by Nash \cite{J.Nash}, has been a foundational concept in the analysis of non-cooperative games. Early works on Nash equilibrium focused primarily on open-loop solutions, such as the conditions for existence and uniqueness explored by Jank and Abou-Kandil \cite{AH} and further extended by Reddy and Zaccour \cite{Reddy} for feedback Nash equilibrium in dynamic games. Over time, scholars like Engwerda \cite{Jacob} investigated the feedback Nash equilibrium properties in scalar LQ differential games, while Moon \cite{Moon} proposed methods to characterize the Nash equilibrium in zero-sum stochastic differential games. Recent studies, such as those by Wang et al. \cite{Wang.el}, delve into non-zero-sum LQ differential games involving asymmetric information structures, while others, like Huang et al. \cite{huang2023linear}, have explored the $\epsilon$-Nash equilibrium in mean-field LQ games. These contributions have paved the way for analyzing more complex systems, particularly in the stochastic context, with a growing focus on closed-loop Nash equilibria.
	
	LQ differential games, a significant class in game theory, have been widely studied in fields such as engineering and economics. The study of deterministic LQ differential games (DLQG) originated with Ho, Bryson, and Baron's \cite{Baron} work on pursuit-evasion games, followed by Schmitendorf's \cite{Schmitendorf} and Bernhard's \cite{Bernhard} explorations of open- and closed-loop strategies. Zhang \cite{zhang2005} later formalized key equivalences for open-loop saddle points, which Delfour \cite{delfour2007} extended to closed-loop cases. The introduction of stochastic LQ games (SLQG) added complexity, as Mou and Yong \cite{Mou} used Hilbert space methods to study these problems, while Sun and Yong \cite{sun2} analyzed open- and closed-loop saddle points. Non-zero-sum LQ games, where Players have independent objectives rather than purely antagonistic goals, further advanced the field. Unlike zero-sum games, where Riccati equations are symmetric, non-zero-sum games feature asymmetric Riccati equations, making the connection between closed-loop and open-loop Nash equilibria more complex \cite{sun2,sun6}. With the inclusion of backward stochastic differential equations (BSDEs,~for short), research in this area continues to grow, offering deeper insights into both zero- and non-zero-sum differential games.

	Discrete-time stochastic systems are crucial for modeling practical problems, particularly in the context of LQ games. Research in this area has made significant progress, with works like Ni et al. \cite{ni2015,ni2016} addressing optimal control solutions and time-inconsistent control in discrete-time mean-field systems. Studies have also extended to stochastic maximum principles, $H_\infty$ filtering, and control problems for both finite and infinite horizons \cite{dong2022maximum,zhang2021robust,zhang2008}. Sun et al. \cite{sunhuiying} contributed to non-cooperative stochastic games, while Zhou \cite{zhou2018} and Gao-Lin \cite{lin1} explored nonzero-sum games and the existence of Nash equilibria in discrete-time LQ systems. These investigations provide essential insights for optimizing control strategies in dynamic decision-making and resource allocation, demonstrating the importance of discrete-time LQ games in both theoretical analysis and practical applications. Nevertheless, in the context of non-zero-sum/zero-sum frameworks, further exploration is needed regarding LQ games with random coefficients. This paper tackles the difficulties of determining the closed-loop Nash equilibrium in discrete-time LQ games with random coefficients, investigating how randomness affects the balanced, interdepended closed-loop strategy.

	This paper examines discrete-time, two-person non-zero-sum stochastic LQ difference games with random coefficients, emphasizing closed-loop solvability and the determination of a closed-loop Nash equilibrium. The results of the relevant parts of our paper is inspired from the \cite{lin1} and \cite{sun1}, which are involved the deterministic coefficients and the former is discussed in discrete time and the latter in continuous time. Our methodology is partially grounded in previous work on discrete-time stochastic LQ optimal control with random coefficients \cite{Wuyiwei2}, the closed-loop solvability of mean-field stochastic LQ control problems \cite{li1}, and mean-field linear-quadratic stochastic differential games \cite{sun5}. Our main contributions are threefold:

(i)\quad\textbf{Novel Riccati Structure:}
It is demonstrated that random coefficients in both the state equation and cost functionals significantly influence the structure of the associated CCREs.
The resulting equation is a fully nonlinear, higher-order backward stochastic difference equation (BS$\Delta$Es, for short),
differing fundamentally from its deterministic LQ game counterparts \cite{lin1}.

(ii) \textbf{Decoupling Framework:}
The discrete-time stochastic systems including nonhomogeneous terms and random coefficients.
It's worth pointing out that the closed-loop Nash equilibrium strategy is derived by decoupling the cross-coupled stochastic Hamiltonian systems,
yielding a structured characterization of equilibrium characteristics.

(iii) \textbf{General Solvability:}
Two Lyapunov-type equations are constructed by the dynamic programming principle approach (DPP, for short),
and the necessary and sufficient condition of closed-loop Nash equilibrium is established,
which contingents upon the regular solvability of two CCREs.

	The remainder of this paper is structured as follows. Section 2 presents the preliminaries. Section 3 gives the key results  related to CCREs in terms of FBS$\Delta$Es and illustrates the state feedback representation of closed-loop Nash equilibrium. Section 4 demonstrate the regular solvability of the associated CCREs. Section 5 concludes the paper.
	
	\section{Preliminaries }\label{sec:2}
	Let $\Omega$ be a fixed sample space  and let $N$ be a positive integer. $\mathcal{T}$ and $\overline{\mathcal{T}}$ represent the  collection of $\{0, \ldots, N-1\}$ and $\{0, \ldots, N\}$ respectively.~Let $(\Omega,\mathcal{F},\mathbb{F},\mathcal{P})$ be a complete filtered probability space with a $\mathcal{P}$-complete filtration $\mathbb{F}:=\left\{\mathcal{F}_{k} \mid k\in \overline{\mathcal{T}}\right\}$ and the one-dimensional scalar noise $\{\omega_k, k \in \overline{\mathcal{T}}\}$ is assumed to be a martingale difference sequence defined on a probability space $(\Omega,\mathcal{F},\mathbb{F},\mathcal{P})$ satisfying:
	\[\mathbb{E}[\omega_k \mid \mathcal{F}_{k-1}] = 0,~~\mathbb{E}[\omega_k^2 \mid \mathcal{F}_{k-1}] = 1.\]
	Furthermore,~let $\mathcal{F}_{k} = \sigma\{\omega_{0}, \omega_{1}, \ldots, \omega_{k}\}$ and $\mathcal{F}_{-1} := \{\emptyset, \Omega\}$.
	
	Throughout this paper, we introduce the notations include that
	\( \mathbb{R}^{n \times m} \) is the space of \( n \times m \) real matrices; \( \mathbb{R}^{n} = \mathbb{R}^{n \times 1} \) and \( \mathbb{R} = \mathbb{R}^{1} \) represent the \( n \)-dimensional real vector space and the set of real numbers, respectively; \( \mathbb{S}^{n} \) is the space of symmetric real matrices of size \( n \times n \); \( \operatorname{tr}(M) \) is the trace of matrix \( M \), defined as the sum of its diagonal elements; \( \langle \cdot, \cdot \rangle \) denotes the inner product in the Euclidean space, where \( \langle M, N \rangle \) is given by \( \operatorname{tr}(M^{\top} N) \); \( \|M\|_{F} \triangleq \sqrt{\operatorname{tr}(M M^{\top})} \) represents the Frobenius norm of matrix \( M \), and \( \mathcal{R}(M) \) represents the range of matrix \( M \); \( M^{\top} \) is the transpose of matrix \( M \); \( M^{\dagger} \) is the Moore-Penrose pseudoinverse of matrix \( M \)
	having the following properties\cite{Penrose}:
	\begin{equation}\nonumber
		\begin{split}
			&M M^{\dagger} M=M, \quad M^{\dagger} M M^{\dagger}=M^{\dagger}, \\&\left(M M^{\dagger}\right)^{\top}=M M^{\dagger}, \quad\left(M^{\dagger} M\right)^{\top}=M^{\dagger} M.
		\end{split}
	\end{equation}
	Furthermore, if \( M \in \mathbb{R}^{n \times m} \) and \( \Psi \in \mathbb{R}^{n \times \ell} \) such that
	\[
	\mathcal{R}(\Psi) \subseteq \mathcal{R}(M),
	\]
	then all the solutions \( \Pi\) to the linear equation
	\[
	M  \Pi=\Psi
	\]
	are given by the following:
	\[
	\Pi=M^{\dagger} \Psi+\left(I-M^{\dagger} M\right) \Gamma, \quad \Gamma \in \mathbb{R}^{m \times \ell}
	\]
	In addition, if \( M^{\top}=M \in \mathbb{S}^{n} \), then
	\[
	M^{\dagger}=\left(M^{\dagger}\right)^{\top}, \quad M M^{\dagger}=M^{\dagger} M ;\] \[\quad \text { and } \quad M \geqslant 0 \Longleftrightarrow M^{\dagger} \geqslant 0.
	\]
	
	For any \( k \in \mathcal{T} \) and Euclidean space \( \mathbb{H} \), we introduce the following spaces used in this paper:
	
	\begin{itemize}
		\item \( L^{2}_{\mathcal{F}_{N-1}}(\Omega; \mathbb{H}) \): the space of \( \mathbb{H} \)-valued \( \mathcal{F}_{N-1} \)-measurable random variables \( \xi \) such that
		\[
		\|\xi\|_{L^{2}_{\mathcal{F}_{N-1}}(\Omega; \mathbb{H})} := \left[\mathbb{E} \|\xi\|_{\mathbb{H}}^2 \right]^{\frac{1}{2}} < \infty.
		\]
		\item \textbf{$L^\infty_{\mathcal{F}_{N-1}}(\Omega; \mathbb{H})$}: the space of all essentially bounded $\mathcal H$-valued, $\mathcal{F}_{N-1}$-measurable random variables.
		\item \( L_{\mathbb{F}}^2\left(\mathcal{T}; \mathbb{H}\right) \): the set of all \( \mathbb{H} \)-valued stochastic processes \( f(\cdot) = \{f_k \mid f_k \text{ is } \mathcal{F}_{k-1} \text{-measurable}, k \in \mathcal{T}\} \) that satisfy
		\[
		\|f(\cdot)\|_{L_{\mathbb{F}}^2\left(\mathcal{T}; \mathbb{H}\right)} := \left[\mathbb{E} \left( \sum_{k=0}^{N-1} \|f_k\|_{\mathbb{H}}^2 \right) \right]^{\frac{1}{2}} < \infty.
		\]
		
		\item \( L_{\mathbb{F}}^{\infty}(\mathcal{T}; \mathbb{H}) \): the space of all \( \mathbb{H} \)-valued essentially bounded stochastic processes \( f(\cdot) = \{f_k \mid f_k \text{ is } \mathcal{F}_{k-1} \text{-measurable}, k \in \mathcal{T}\} \).
	\end{itemize}
	
	Consider the following linear controlled  stochastic difference equation on finite horizon $\mathcal{T}$:
	\begin{equation}\label{eq:1.1}
		\left\{\begin{array}{ll}
			\begin{split}
				x_{k+1}=&A_{k}x_{k}+B_{k}u_{k}+C_{k}v_{k}+b_{k}\\&+\left(D_{k}x_{k}+E_{k}u_{k}+F_{k}v_{k}+\sigma_{k}\right)\omega_{k}, \\
				x_{0} =&\xi \in\mathbb{R}^{n}~~~~k \in \mathcal{T},
			\end{split}
		\end{array}\right.
	\end{equation}
	where the coefficients $ A_{k}, B_{k}, C_{k}, D_{k}, E_{k}, F_{k}$ are \(\mathcal{F}_{k-1}\)-measurable random coefficients, and \( b_{k}, \sigma_{k} \) are \(\mathcal{F}_{k-1}\)-measurable vector-valued processes.
	
	The admissible control spaces for Players 1 and 2 are defined as follows:
	
	\begin{defn}[Control Spaces]
		\label{def:control_spaces}
		\begin{subequations}
			\begin{align}
				\mathcal{U} &\equiv \mathcal{U}(\mathcal{T}) \nonumber \\
				&\equiv \Big\{ \mathbf{u} = (u_0, u_1, \ldots, u_{N-1}) \,\Big|\, \forall k \in \mathcal{T},\ u_k: \Omega \to \mathbb{R}^m \text{ is } \nonumber \\
				&\quad \begin{aligned}
					&\mathcal{F}_{k-1}\text{-measurable};~\mathbb{E} \bigg[ \sum_{k=0}^{N-1} \|u_k\|^2 \bigg] < \infty
				\end{aligned}\Big\},
				\label{eq:control_space_U} \\
				\mathcal{V} &\equiv \mathcal{V}(\mathcal{T}) \nonumber \\
				&\equiv \Big\{ \mathbf{v} = (v_0, v_1, \ldots, v_{N-1}) \,\Big|\, \forall k \in \mathcal{T},\ v_k: \Omega \to \mathbb{R}^l \text{ is }\nonumber \\
				&\quad \begin{aligned}
					&\mathcal{F}_{k-1}\text{-measurable}; ~\mathbb{E} \bigg[ \sum_{k=0}^{N-1} \|v_k\|^2 \bigg] < \infty
				\end{aligned} \Big\}.
				\label{eq:control_space_V}
			\end{align}
		\end{subequations}
	\end{defn}
	
	To evaluate the control performance of $\mathbf{u}$ and $\mathbf{v}$, we define the cost functionals for both Players:
	\begin{equation}\label{eq:1.4}
		\begin{split}
			\mathcal{J}_1&(0, \xi; \mathbf{u}, \mathbf{v}) \\=& \frac{1}{2} \mathbb{E}\bigg\{\langle G_{N} x_{N}, x_{N}\rangle + 2\langle g_{N}, x_{N}\rangle + \sum_{k=0}^{N-1} \bigg[\langle Q_{k}x_{k}, x_{k}\rangle\\
			&+ 2\langle L_{k}^{\top} u_{k}, x_{k}\rangle + \langle R_{k}u_{k}, u_{k}\rangle  + 2\langle q_{k}, x_{k}\rangle + 2\langle \rho_{k}, u_{k}\rangle\bigg]\bigg\},
		\end{split}
	\end{equation}
	for Player 1 and
	\begin{equation}\label{eq:1.5}
		\begin{split}
			\mathcal{J}_2&(0, \xi; \mathbf{u}, \mathbf{v}) \\=& \frac{1}{2} \mathbb{E}\bigg\{\langle H_{N} x_{N}, x_{N}\rangle + 2\langle h_{N}, x_{N}\rangle +\sum_{k=0}^{N-1}  \bigg[\langle P_{k}x_{k}, x_{k}\rangle \\
			&+ 2\langle M_{k}^{\top} v_{k}, x_{k}\rangle + \langle S_{k}v_{k}, v_{k}\rangle + 2\langle p_{k}, x_{k}\rangle + 2\langle \theta_{k}, v_{k}\rangle\bigg]\bigg\},
		\end{split}
	\end{equation}
	for Player 2.~Here,~\(G_{N}\) and \(H_{N}\) are \(\mathcal{F}_{N-1}\)-measurable random symmetric real matrices; The random weighting coefficients \( Q_{k}, P_{k}, L_{k}, M_{k}, R_{k}, S_{k} \) are \(\mathcal{F}_{k-1}\)-measurable matrix-valued functions with \(Q^{\top}_{k} = Q_{k}\), \(R^{\top}_{k} = R_{k}\), \(P^{\top}_{k} = P_{k}\), \(S^{\top}_{k} = S_{k}\);  \(g_{N},~h_{N}\) \(\in L^{2}_{\mathcal{F}_{N-1}}(\Omega; \mathbb{R}^{n})\); and \(q_{k}, \rho_{k}, p_{k}, \theta_{k}\) are vector-valued \(\mathcal{F}_{k-1}\)-measurable processes.

	To ensure the well-posedness of the closed-loop system and the admissibility of feedback controls, we introduce the following key functional spaces with essential boundedness conditions:
	
	\[
	\begin{cases}
		\begin{aligned}
			\mathcal{Q}^1(\mathcal{T}) &:= L_{\mathbb{F}}^{\infty}(\mathcal{T}; \mathbb{R}^{m \times n}) \\
			\mathcal{Q}^2(\mathcal{T}) &:= L_{\mathbb{F}}^{\infty}(\mathcal{T}; \mathbb{R}^{l \times n}) \\
			\mathscr{M}^1(\mathcal{T}) &:= \mathcal{Q}^1(\mathcal{T}) \times \mathcal{U} \\
			\mathscr{M}^2(\mathcal{T}) &:= \mathcal{Q}^2(\mathcal{T}) \times \mathcal{V} \\
			\mathscr{M}(\mathcal{T}) &:= \mathscr{M}^1(\mathcal{T}) \times \mathscr{M}^2(\mathcal{T})
		\end{aligned}
	\end{cases}\]
	
	Next, we impose the following Assumptions about the state equation $\eqref{eq:1.1}$, cost functionals $\eqref{eq:1.4}$ and $\eqref{eq:1.5}$:
	\begin{ass}\label{ass:2.1}
		The coefficients \( A_{k}, D_{k}: \Omega \rightarrow \mathbb{R}^{n\times n} ,\) \(B_{k}, E_{k}: \Omega \rightarrow \mathbb{R}^{n\times m} \), \( C_{k}, F_{k}: \Omega \rightarrow \mathbb{R}^{n\times l} \) are uniformly bounded \(\mathcal{F}_{k-1}\)-measurable mappings, and \( b., \sigma. \in L_{\mathbb{F}}^2(\mathcal{T}; \mathbb{R}^{n}) \).
	\end{ass}

	\begin{ass}\label{ass:2.2}
		The weighting coefficients \( Q_{k}, P_{k}: \Omega \rightarrow \mathbb{R}^{n\times n} \), \( L_{k}: \Omega \rightarrow \mathbb{R}^{m\times n} \), \( M_{k}: \Omega \rightarrow \mathbb{R}^{l\times n} \), \( R_{k}: \Omega \rightarrow \mathbb{R}^{m\times m} \), \( S_{k}: \Omega \rightarrow \mathbb{R}^{l\times l} \) are uniformly bounded \(\mathcal{F}_{k-1}\)-measurable mappings; \( G_{N}, H_{N}\in L^\infty_{\mathcal{F}_{N-1}}(\Omega; \mathbb{R}^{n}) \) ; \( g_{N}, h_{N} \in L_{\mathcal{F}_{N-1}}^2(\Omega; \mathbb{R}^{n}) \), \( q., p. \in L_{\mathbb{F}}^2(\mathcal{T}; \mathbb{R}^{n}) \), \( \rho. \in L_{\mathbb{F}}^2(\mathcal{T}; \mathbb{R}^{m}) \), \( \theta. \in L_{\mathbb{F}}^2(\mathcal{T}; \mathbb{R}^{l}) \).
	\end{ass}
	\begin{ass}\label{ass:2.3}
		For all $ k \in \mathcal{T} $, there exists a constant $\delta > 0$ such that
		\begin{equation}\label{eq:20.6}
			\left\{
			\begin{aligned}
				&G_{N},~H_{N} \succeq 0,~~ R_{k},~S_{k} \succeq \delta I, \\
				&Q_{k} - L^{\top}_{k} R^{-1}_{k} L_{k} \succeq 0,~~ \\
				&P_{k} - M^{\top}_{k} S^{-1}_{k} M_{k} \succeq 0, \quad k \in \mathcal{T}.
			\end{aligned}
			\right.
		\end{equation}
	\end{ass}
	Given Assumption \ref{ass:2.1}, by \cite{lin1}, for any initial state \( \xi \in \mathbb{R}^{n} \), \( \mathbf{u} \in \mathcal{U} \), and \( \mathbf{v} \in \mathcal{V} \), the state equation \eqref{eq:1.1} has a unique solution denoted as $ x_{(\cdot)} \equiv \{ x^{(\xi, \mathbf{u}, \mathbf{v})}_k \}_{k=0}^N \in L_{\mathbb{F}}^2(\overline{\mathcal{T}}; \mathbb{R}^{n}) $. Furthermore, under Assumptions \ref{ass:2.1} and \ref{ass:2.2}, for any initial state \( \xi \in \mathbb{R}^{n} \), \( \mathbf{u} \in \mathcal{U} \), and \( \mathbf{v} \in \mathcal{V} \), it is straightforward to verify that the cost functionals \eqref{eq:1.4} and \eqref{eq:1.5} are well-defined. Therefore, the discrete-time two-person non-zero-sum stochastic difference games with random coefficients is well-defined.~Next, we will illustrate two nonzero-sum Nash equilibrium problems formally by mathematical terms.
	\begin{pro} [\textbf{DSDG-OL}]\label{pro:1.1}
		~~\\
		Consider the discrete-time stochastic difference game \eqref{eq:1.1} with any initial state \( \xi \in \mathbb{R}^{n} \), and admissible control sets \( \mathcal{U}, \mathcal{V} \) defined in \eqref {eq:control_space_U} and \eqref {eq:control_space_V}.~The objective is to find an admissible control strategy pair \( \left(\mathbf{u}^{*}, \mathbf{v}^{*}\right) \in \mathcal{U} \times \mathcal{V} \) such that:
		\begin{equation}\label{eq:1.6}
			\left\{
			\begin{array}{ll}
				\mathbf{u}^{*} \in \underset{\mathbf{u} \in \mathcal{U}}{\arg\min}  \mathcal{J}_{1}\left(0, \xi ; \mathbf{u}, \mathbf{v}^{*}\right),
				\\
				\mathbf{v}^{*} \in \underset{\mathbf{v} \in \mathcal{V}}{\arg\min} \mathcal{J}_{2}\left(0, \xi ; \mathbf{u}^{*}, \mathbf{v}\right) .
			\end{array}
			\right.
		\end{equation}
	\end{pro}
	Any \(  (\mathbf{u}^{*}, \mathbf{v}^{*} )\)  \( \in \mathcal{U} \times \mathcal{V} \) satisfying \eqref{eq:1.6} is called an \textbf{open-loop Nash equilibrium} of Problem (DSDG-OL), which is solved in full in\cite{Wuyiwei}.~Correspondingly, the formal statement for closed-loop Nash equilibrium is in the manner outlined below.
	
	\begin{pro} [\textbf{DSDG-CL}]\label{pro:1.2}
		~~\\
		Consider the discrete-time stochastic difference game \eqref{eq:1.1} with any initial state \( \xi \in \mathbb{R}^{n} \).~The goal is to find a quadruple strategy	
		$$\left(\Pi^{1 *}, \Sigma^{1 *}; \Pi^{2 *}, \Sigma^{2 *}\right) \in \mathscr{M}(\mathcal{T})$$
		such that the following conditions hold:
		\begin{equation}\label{eq:2.8}
			\left\{
			\begin{array}{ll}
				&V_1(0, \xi) \\&\triangleq \mathcal{J}_1(0, \xi, \Pi^{1*}\mathbf{x}^{\ast}+\Sigma^{1\ast}, \Pi^{2*}\mathbf{x}^{\ast}+\Sigma^{2\ast}) \\&= \inf_{\big(\Pi^{1},\Sigma^{1}\big) \in \mathscr{M}^{1}(\mathcal{T})} \mathcal{J}_1(0, \xi, \Pi^{1}\mathbf{x}_{1}+\Sigma^{1}, \Pi^{2*}\mathbf{x}_{1}+\Sigma^{2\ast}),
				\\
				&V_2(0, \xi) \\&\triangleq \mathcal{J}_2(0, \xi, \Pi^{1*}\mathbf{x}^{\ast}+\Sigma^{1\ast}, \Pi^{2*}\mathbf{x}^{\ast}+\Sigma^{2\ast}) \\&= \inf_{\big(\Pi^{2},\Sigma^{2}\big) \in \mathscr{M}^{2}(\mathcal{T})} \mathcal{J}_2(0, \xi, \Pi^{1*}\mathbf{x}_{2}+\Sigma^{1\ast}, \Pi^{2}\mathbf{x}_{2}+\Sigma^{2}),
			\end{array}
			\right.
		\end{equation}
		where $\mathbf{x}^* = X(\Pi^{1*},\Sigma^{1*};\Pi^{2*},\Sigma^{2*},\xi)$ denotes the equilibrium state trajectory generated by the closed-loop control law
		\begin{equation}\label{eq:2.9}
			u_k^* = \Pi_k^{1*}x_k^* + \Sigma_k^{1*}, \quad v_k^* = \Pi_k^{2*}x_k^* + \Sigma_k^{2*}, \quad k \in \mathcal{T}.
		\end{equation}
		The trajectory $\mathbf{x}_1 = X(\Pi^1,\Sigma^1;\Pi^{2*},\Sigma^{2*},\xi)$, corresponding to Player 1's cost functional, is generated by
		\begin{equation}\label{eq:2.100}
			u_k = \Pi_k^1x_k + \Sigma_k^1, \quad v_k = \Pi_k^{2*}x_k + \Sigma_k^{2*}, \quad k \in \mathcal{T},
		\end{equation}
		while $\mathbf{x}_2 = X(\Pi^{1*},\Sigma^{1*};\Pi^2,\Sigma^2,\xi)$, corresponding to Player 2's cost functional, is generated by
		\begin{equation}\label{eq:2.110}
			u_k = \Pi_k^{1*}x_k + \Sigma_k^{1*}, \quad v_k = \Pi_k^2x_k + \Sigma_k^2, \quad k \in \mathcal{T}.
		\end{equation}

		A quadruple $\big(\Pi^{1\ast},\Sigma^{1\ast};\Pi^{2\ast},\Sigma^{2\ast}\big)$ satisfying \eqref{eq:2.8} is called a \textbf{closed-loop Nash equilibrium} of Problem (DSDG-CL) on finite horizon $\mathcal{T}$ if for any $\xi \in \mathbb{R}^{n}$ and any quadruple $\big(\Pi^{1},\Sigma^{1};\Pi^{2},\Sigma^{2}\big)\in\mathscr{M}( \mathcal{T})$.~Furthermore,~Problem (DSDG-CL) is denoted as (DSDG-CL)$^{0}$ when Problem (DSDG-CL) omitting the nonhomogeneous terms $b_{},\sigma_{},g_{N},q_{},\rho_{},h_{N},p_{},\theta_{}$. \end{pro}

	\begin{rmk}
		~~~~~~\\
		(i) If a closed-loop Nash equilibrium $\big(\Pi^{1*},\Sigma^{1*};\Pi^{2*},\Sigma^{2*}\big)$ (uniquely) exists on $\mathcal{T}$, Problem 1.2 is said to be (uniquely) closed-loop solvable on $\mathcal{T}$.
		~~~~~~\\
		(ii) An open-loop Nash equilibrium \(  (\mathbf{u}^{*}, \mathbf{v}^{*} )\) is allowed to depend on the initial state \(\xi\), whereas the closed-loop Nash equilibrium  $\big(\Pi^{1*},\Sigma^{1*};\Pi^{2*},\Sigma^{2*}\big)$ is required to be independent of the initial state $\xi$.
		~~~~~~\\
		(iii) Let control law \eqref{eq:2.9} holds, it's clear that the closed-loop Nash equilibrium strategy pair \( \big(\mathbf{u}^{*},\mathbf{v}^{*}\big) \equiv \big(\Pi^{1*} \mathbf{x}^{*}+ \Sigma^{1*} , \Pi^{2*} \mathbf{x}^{*} + \Sigma^{2*} \big)\) induced by the closed-loop Nash equilibrium \(\left(\Pi^{1*}, \Sigma^{1*}; \Pi^{2*}, \Sigma^{2*}\right)\) is simultaneously an open-loop Nash equilibrium of Problem (DSDG-OL) for the initial value \(\xi\). Hence, closed-loop solvability implies open-loop solvability but the converse is not necessarily true.
	\end{rmk}
	
	From the result of Sun and Yong\cite{sun2}, there have an important proposition providing some equivalent definitions of closed-loop Nash equilibrium.~We present below a similar proposition that characterizes the Closed-Loop Nash Equilibrium property:
	\begin{prop}\label{prop:nash_equilibrium}(Characterization of Closed-Loop Nash Equilibrium) Let Assumptions \ref{ass:2.1}-\ref{ass:2.2} hold. Given a quadruple strategy of feedback control law
		$$\left(\Pi^{1 *}, \Sigma^{1 *}; \Pi^{2 *}, \Sigma^{2 *}\right) \in \mathscr{M}(\mathcal{T}),$$
		and let the corresponding state trajectory of system $\eqref{eq:1.1}$ under this  control law be denoted by $\mathbf{x}^*$, then the following statements are equivalent:
		~~~~\\
		(i) \(\left(\Pi^{1*}, \Sigma^{1*}; \Pi^{2*}, \Sigma^{2*}, \xi\right)\) is a closed-loop Nash equilibrium of Problem (DSDG-CL) on finite horizon $\mathcal{T}$.
		~~~\\
		(ii) For any $\xi\in\mathbb{R}^{n}$, \( \Sigma^1\in\mathcal{U}\) and the corresponding trajectory $\mathbf{x}$,
		\begin{equation}\label{eq:2.10}
			\begin{split}
				\mathcal{J}&_1\left(0, \xi; \Pi^{1*} \mathbf{x}^{*} + \Sigma^{1*}, \Pi^{2*} \mathbf{x}^{*}  + \Sigma^{2*}\right) \\\leq \mathcal{J}&_1\left(0, \xi ; \Pi^{1*} \mathbf{x}_1  + \Sigma^1, \Pi^{2*} \mathbf{x}_1 + \Sigma^{2*}\right),
			\end{split}
		\end{equation}
		
		For any $\xi\in\mathbb{R}^{n}$, \( \Sigma^2 \in \mathcal{V}\) and the corresponding trajectory $\mathbf{x}$,	    	\begin{equation}\label{eq:2.11}
			\begin{split}
				\mathcal{J}&_2\left(0, \xi; \Pi^{1*} \mathbf{x}^{*} + \Sigma^{1*}, \Pi^{2*} \mathbf{x}^{*}  + \Sigma^{2*}\right) \\\leq \mathcal{J}&_2\left(0, \xi ; \Pi^{1*} \mathbf{x}_2 + \Sigma^{1*}, \Pi^{2*} \mathbf{x}_2 + \Sigma^2\right),
			\end{split}
		\end{equation}
		where $\mathbf{x}^*$ denotes the equilibrium state trajectory generated by \eqref{eq:2.9}, the trajectory $\mathbf{x}_1$ in \eqref{eq:2.10} is generated by \eqref{eq:2.100}  and the trajectory $\mathbf{x}_2$ in \eqref{eq:2.11} is generated by \eqref{eq:2.110}.
		
		(iii) For any $\xi\in\mathbb{R}^{n}$,~any admissible control pair \((\mathbf{u}, \mathbf{v}) \in \mathcal{U} \times \mathcal{V}\), and the corresponding trajectory $\mathbf{x}_1,\mathbf{x}_2$,
		\begin{equation}\label{eq:2.12}
			\mathcal{J}_1(0, \xi ; \Pi^{1*} \mathbf{x}^{*} + \Sigma^{1*}, \Pi^{2*} \mathbf{x}^{*} + \Sigma^{2*}) \leq \mathcal{J}_1(0, \xi ; \mathbf{u}, \Pi^{2*} \mathbf{x}_1 + \Sigma^{2*}),
		\end{equation}
		\begin{equation}\label{eq:2.13}
			\mathcal{J}_2(0, \xi ; \Pi^{1*} \mathbf{x}^{*} + \Sigma^{1*}, \Pi^{2*} \mathbf{x}^{*} + \Sigma^{2*}) \leq \mathcal{J}_2(0, \xi ; \Pi^{1*} \mathbf{x}_2 + \Sigma^{1*}, \mathbf{v}),
		\end{equation}
		where $\mathbf{x}$ is the equilibrium state trajectory generated by the closed-loop law \eqref{eq:2.9}, $\mathbf{x}_1$ and $\mathbf{x}_2$ are generated by \eqref{eq:2.100} and \eqref{eq:2.110} respectively.
	\end{prop}

	\begin{proof}
		\textbf{(i) $\Rightarrow$ (ii):}
		Direct from Definition \ref{pro:1.2} by fixing $\Pi^2 = \Pi^{2*}$ in Player 1's optimization and $\Pi^1 = \Pi^{1*}$ in Player 2's optimization.
		
		\textbf{(ii) $\Rightarrow$ (iii):}
		Fix $\xi \in \mathbb{R}^n$ and $(\mathbf{u}, \mathbf{v}) \in \mathcal{U} \times \mathcal{V}$. Consider the closed-loop system:
		\begin{equation*}\label{eq:modified_system}
			\left\{\begin{aligned}
				x_{k+1} &= (A_k + B_k \Pi_k^{1*} + C_k \Pi_k^{2*}) x_k + B_k (u_k - \Pi_k^{1*} x_k) \\
				&+ C_k \Sigma_k^{2*} + b_k + \left[(D_k + E_k \Pi_k^{1*} + F_k \Pi_k^{2*}) x_k \right. \\
				&+ \left. E_k (u_k - \Pi_k^{1*} x_k) + F_k \Sigma_k^{2*} + \sigma_k \right] \omega_k \\
				x_0 &= \xi
			\end{aligned}\right.
		\end{equation*}
		
		Define $\Sigma_k^1 \triangleq u_k - \Pi_k^{1*} x_k$. This is $\mathcal{F}_{k-1}$-measurable because:
		\begin{itemize}
			\item $u_k$ is $\mathcal{F}_{k-1}$-measurable (since $\mathbf{u} \in \mathcal{U}$)
			\item $\Pi_k^{1*}$ is $\mathcal{F}_{k-1}$-measurable (since $\Pi^{1*} \in \mathcal{Q}^1(\mathcal{T}) \subset L^\infty_{\mathbb{F}}$)
			\item $x_k$ is $\mathcal{F}_{k-1}$-measurable (adapted state process)
			\item $\Sigma^1 \in \mathcal{U}$ as $|u_k|^2 \leq 2(|\Pi_k^{1*} x_k|^2 + |\Sigma_k^1|^2)$ and $\|\Pi^{1*}\|_\infty < \infty$ by $L^\infty$ assumption
		\end{itemize}
		
		The trajectory $\mathbf{x}$ coincides with $\mathbf{x}_1$ from \eqref{eq:2.100} under controls:
		$$
		u_k = \Pi_k^{1*} x_k + \Sigma_k^1, \quad v_k = \Pi_k^{2*} x_k + \Sigma_k^{2*}
		$$
		Applying \eqref{eq:2.10}:
		\begin{align*}
			\mathcal{J}_1&\left(0, \xi; \Pi^{1*} \mathbf{x}^{*} + \Sigma^{1*}, \Pi^{2*} \mathbf{x}^{*}  + \Sigma^{2*}\right) \\
			&\leq \mathcal{J}_1\left(0, \xi; \Pi^{1*} \mathbf{x} + \Sigma^1, \Pi^{2*} \mathbf{x} + \Sigma^{2*}\right) \\
			&= \mathcal{J}_1\left(0, \xi; \mathbf{u}, \Pi^{2*} \mathbf{x} + \Sigma^{2*}\right)
		\end{align*}
		This proves \eqref{eq:2.12}. The Player 2 case \eqref{eq:2.13} is symmetric.
		
		\textbf{(iii) $\Rightarrow$ (i):}
		Fix $(\Pi^2, \Sigma^2) = (\Pi^{2*}, \Sigma^{2*})$ and take any $(\Pi^1, \Sigma^1) \in \mathscr{M}^1(\mathcal{T})$. Define $\mathbf{u} \triangleq \Pi^1\mathbf{x}_1 + \Sigma^1 \in \mathcal{U}$ (since $\Pi^1 \in L^\infty_\mathbb{F}$ and $\Sigma^1 \in \mathcal{U}$). Let $\mathbf{x}_1$ be generated by:
		$$
		\left\{\begin{aligned}
			u_k &= \Pi_k^1 x_k^1 + \Sigma_k^1 \\
			v_k &= \Pi_k^{2*} x_k^1 + \Sigma_k^{2*}
		\end{aligned}\right.
		$$

		Apply \eqref{eq:2.12} with $\mathbf{x}_1$ and $\mathbf{u} = \Pi^1\mathbf{x}_1 + \Sigma^1$:
		\begin{align*}
			&\mathcal{J}_1(0,\xi; \Pi^{1*}\mathbf{x}^* + \Sigma^{1*}, \Pi^{2*}\mathbf{x}^* + \Sigma^{2*})
			\\&\leq \mathcal{J}_1(0,\xi; \mathbf{u}, \Pi^{2*}\mathbf{x}_1 + \Sigma^{2*}) \\
			&= \mathcal{J}_1(0,\xi; \Pi^1\mathbf{x}_1 + \Sigma^1, \Pi^{2*}\mathbf{x}_1 + \Sigma^{2*})
		\end{align*}
		This proves Player 1's condition. Player 2's proof is symmetric.
	\end{proof}

	\section{FBS$\Delta$Es and CCREs in Closed-Loop Form}\label{sec:3}
	\normalsize{
		In this section, the closed-loop Nash equilibrium is studied in terms of FBS$\Delta$Es, which is a well-known stochastic maximum principle (SMP, for short) approach in solving difference games, as it combines the state equations, cost functionals, and optimality conditions into a unified framework. To establish the closed-loop Nash equilibrium solution, we first reformulate the original stochastic Hamiltonian system from \cite{Wuyiwei}, specifically in section 4.1, in terms of feedback strategies. The transition from open-loop to closed-loop control structure necessitates the incorporation of the closed-loop feedback control law \eqref{eq:2.9} into the Hamiltonian framework. Afterthat, we will present two sets of CCREs to decouple the Hamiltonian systems with feedback strategies.
		
		The Hamiltonian functions \( \mathcal{H}_1 \) and \( \mathcal{H}_2 \) for Player 1 and Player 2, respectively, are defined by combining the state dynamics in \eqref{eq:1.1} with the cost functionals in \eqref{eq:1.4} and \eqref{eq:1.5}, and are augmented by the dual variables \( y_{1}=\{y_{1,k}\}_{k \in \mathcal{T}} \in L_{\mathbb{F}}^2(\overline{\mathcal{T}}; \mathbb{R}^n) \) and \( y_{2}=\{y_{2,k}\}_{k \in \mathcal{T}} \in L_{\mathbb{F}}^2(\overline{\mathcal{T}}; \mathbb{R}^n) \), respectively. Specifically, the Hamiltonian function for Player 1 is given by:
		\[
		\begin{aligned}
			&\mathcal{H}_1(k, x_k, u_k, v_k, y_{1,k+1}) \\=& \frac{1}{2} \big( \langle Q_k x_k, x_k \rangle + 2 \langle L_k^\top u_k, x_k \rangle + \langle R_k u_k, u_k \rangle + 2 \langle q_k, x_k \rangle\\
			&  + 2 \langle \rho_k, u_k \rangle \big) + \langle \mathbb{E}\big[y_{1,k+1}\mid \mathcal{F}_{k-1} \big], A_k x_k + B_k u_k + C_k v_k \\
			&+ b_k  + (D_k x_k + E_k u_k + F_k v_k + \sigma_k) \omega_k \rangle ,
		\end{aligned}
		\]
		and for Player 2, the Hamiltonian function is defined as:
		\[
		\begin{aligned}
			&\mathcal{H}_2(k, x_k, u_k, v_k, y_{2,k+1}) \\=& \frac{1}{2} \big( \langle P_k x_k, x_k \rangle + 2 \langle M_k^\top v_k, x_k \rangle + \langle S_k v_k, v_k \rangle + 2 \langle p_k, x_k \rangle \\
			& + 2 \langle \theta_k, v_k \rangle \big) + \langle\mathbb{E}\big[ y_{2,k+1}\mid \mathcal{F}_{k-1} \big], A_k x_k + B_k u_k + C_k v_k \\
			&+ b_k  + (D_k x_k + E_k u_k + F_k v_k + \sigma_k) \omega_k \rangle .
		\end{aligned}
		\]
		Under feedback strategies of the form $v_k = \Pi_k^2x_k + \Sigma_k^2$, the Hamiltonian function $\mathcal{H}_1(k, x_k, u_k, v_k, y_{1,k+1})$ transforms into the augmented form by the feedback strategies of the form $v_{k}=\Pi^{2}_{k}x_{k}+\Sigma^{2}_{k}$:
		\begin{equation}\label{eq:3.1}
			\begin{split}
				&\mathcal{H}_1(k, x_k, u_k, \Pi^{2}_{k}x_{k}+\Sigma^{2}_{k}, y_{1,k+1})
				\\= & \frac{1}{2} \big( \langle Q_k x_k, x_k \rangle + 2 \langle L_k^\top u_k, x_k \rangle + \langle R_k u_k, u_k \rangle  + 2 \langle q_k, x_k \rangle \\
				&+ 2 \langle \rho_k, u_k \rangle \big)  + \langle \mathbb{E}\big[y_{1,k+1}\mid \mathcal{F}_{k-1}\big], A_k x_k + B_k u_k + C_k\\
				&\times (\Pi^{2}_{k}x_{k}+\Sigma^{2}_{k})  + b_k + (D_k x_k + E_k u_k + F_k (\Pi^{2}_{k}x_{k}+\Sigma^{2}_{k}) \\&+ \sigma_k) \omega_k \rangle.
			\end{split}
		\end{equation}
		And the Hamiltonian function  $\mathcal{H}_2(k, x_k, u_k, v_k, y_{2,k+1})$ transforms into the augmented form by the feedback strategies of the form $u_{k}=\Pi^{1}_{k}x_{k}+\Sigma^{1}_{k}$:
		\begin{equation}\label{eq:3.2}
			\begin{split}
				&\mathcal{H}_2(k, x_k, \Pi^{1}_{k}x_{k}+\Sigma^{1}_{k}, v_k, y_{2,k+1})
				\\= & \frac{1}{2} \big( \langle P_k x_k, x_k \rangle + 2 \langle M_k^\top v_k, x_k \rangle + \langle S_k v_k, v_k \rangle  + 2 \langle p_k, x_k \rangle \\
				&+ 2 \langle \theta_k, v_k \rangle \big)  + \langle \mathbb{E}\big[y_{2,k+1}\mid \mathcal{F}_{k-1}\big], A_k x_k + B_k (\Pi^{1}_{k}x_{k}+\Sigma^{1}_{k}) \\
				&+ C_k v_k + b_k  + (D_k x_k + E_k (\Pi^{1}_{k}x_{k}+\Sigma^{1}_{k}) + F_k v_k + \sigma_k) \omega_k \rangle.
			\end{split}
		\end{equation}
		The adjoint processes, denoted as \( \{y_{1,k}\}_{k \in \mathcal{T}} \) for Player~1, corresponding to the admissible triple \((\mathbf{x}, \mathbf{u}, \Pi^{2}\mathbf{x}+\Sigma^{2})\), and \( \{y_{2,k}\}_{k \in \mathcal{T}} \) for Player~2, corresponding to the admissible triple \((\mathbf{x}, \Pi^{1}\mathbf{x}+\Sigma^{1}, \mathbf{v})\), are obtained by differentiating their respective Hamiltonians with respect to \(x_k\).~For Player 1, the adjoint process \( y_{1}=\{y_{1,k}\}_{k \in \mathcal{T}} \) is given by:
		\begin{equation}\label{eq:3.3}
			\begin{split}
				y_{1,k} =& \frac{\partial \mathcal{H}_1(k, x_k, u_k, \Pi^{2}_{k}x_{k}+\Sigma^{2}_{k}, y_{1,k+1}) }{\partial x_k} \\
				=& \big(A_{k}+C_{k}\Pi^{2}_{k}\big)^{\top}\mathbb{E}[y_{1, k+1} \mid \mathcal{F}_{k-1}] + \big(D_{k}+F_{k}\Pi^{2}_{k}\big)^{\top}\\&\times\mathbb{E}[y_{1, k+1}\omega_{k} \mid \mathcal{F}_{k-1}] + Q_{k}x_{k} + L_{k}^{\top}u_{k}+q_{k}.
			\end{split}
		\end{equation}
		Similarly, for Player 2, the adjoint process \( y_{2}=\{y_{2,k}\}_{k \in \mathcal{T}} \)  is given by:
		\begin{equation}\label{eq:3.4}
			\begin{split}
				y_{2,k} =& \frac{\partial \mathcal{H}_2(k, x_k, \Pi^{1}_{k}x_{k}+\Sigma^{1}_{k}, v_k, y_{2,k+1})}{\partial x_k} \\
				=& \big(A_{k}+B_{k}\Pi^{1}_{k}\big)^{\top}\mathbb{E}[y_{2, k+1} \mid \mathcal{F}_{k-1}] + \big(D_{k}+E_{k}\Pi^{1}_{k}\big)^{\top}\\&\times\mathbb{E}[y_{2, k+1}\omega_{k} \mid \mathcal{F}_{k-1}] + P_{k}x_{k} + M_{k}^{\top}v_{k}+p_{k}.
			\end{split}
		\end{equation}
		Furthermore,~we can derive the stationary
		conditions that the Nash equilibrium must satisfy in the following

		\begin{thm}\label{thm:3.1A}
			Let Assumptions \ref{ass:2.1}--\ref{ass:2.3} hold. Then, quadruple strategy
			$$
			\left(\Pi^{1*}, \Sigma^{1*}; \Pi^{2*}, \Sigma^{2*}\right) \in \mathscr{M}(\mathcal{T})
			$$
			is a closed-loop Nash equilibrium of Problem~\ref{pro:1.2} on the finite horizon $\mathcal{T}$, if and only if the following two stationarity conditions hold simultaneously:
			
			\begin{enumerate}
				\item For Player 1, given Player 2's equilibrium strategy pair $\left(\Pi^{2*}, \Sigma^{2*}\right) \in \mathscr{M}^{2}(\mathcal{T}) $, the following stationarity condition holds for all $k \in \mathcal{T}$:
				\begin{equation}\label{eq:3.5}
					\begin{split}
						&\frac{\partial \mathcal{H}_1(k, x^{*}_k, u^{*}_k, \Pi^{2*}_{k}x^*_{k}+\Sigma^{2*}_{k}, y_{1,k+1})}{\partial u^{*}_k} \\
						=~& B_k^\top \mathbb{E}[y_{1,k+1} \mid \mathcal{F}_{k-1}] + E_k^\top \mathbb{E}[y_{1,k+1} \omega_k \mid \mathcal{F}_{k-1}] \\&+ L_k x^{*}_k + R_k u^{*}_k + \rho_k = 0,
					\end{split}
				\end{equation}
				where adjoint process \(y_1 = \{y_{1,k}\}_{k\in \mathcal{T}}\)  satisfies:
				\begin{equation}\label{eq:3.6}
					\begin{split}
						y_{1,k} = &(A_k + C_k \Pi_k^{2*})^\top \mathbb{E}[y_{1,k+1} \mid \mathcal{F}_{k-1}] \\
						&+ (D_k + F_k \Pi_k^{2*})^\top \mathbb{E}[y_{1,k+1} \omega_k \mid \mathcal{F}_{k-1}] \\
						&+ Q_k x_k^* + L_k^\top u_k^* + q_k,
					\end{split}
				\end{equation}
				and $u^{*}_{k}=\Pi^{1*}_{k}x^{*}_{k}+\Sigma^{1*}_{k}$ in both \eqref{eq:3.5} and \eqref{eq:3.6}.
				\item For Player 2, given Player 1's equilibrium strategy pair $\left(\Pi^{1*}, \Sigma^{1*}\right) \in \mathscr{M}^{1}(\mathcal{T}) $, the following stationarity condition holds for all $k \in \mathcal{T}$:
				\begin{equation}\label{eq:3.7}
					\begin{split}
						&\frac{\partial \mathcal{H}_2(k, x^{*}_k, \Pi^{1*}_{k}x^*_{k}+\Sigma^{1*}_{k}, v^{*}_k, y_{2,k+1})}{\partial v^{*}_k} \\
						=~& C_k^\top \mathbb{E}[y_{2,k+1} \mid \mathcal{F}_{k-1}] + F_k^\top \mathbb{E}[y_{2,k+1} \omega_k \mid \mathcal{F}_{k-1}] \\&+ M_k x^{*}_k + S_k v^{*}_k + \theta_k = 0,
					\end{split}
				\end{equation}
				where adjoint process \( y_2 = \{y_{2,k}\}_{k\in \mathcal{T}}\)  satisfies:
				\begin{equation}\label{eq:3.8}
					\begin{split}
						y_{2,k} = &(A_k + B_k \Pi_k^{1*})^\top \mathbb{E}[y_{2,k+1} \mid \mathcal{F}_{k-1}] \\
						&+ (D_k + E_k \Pi_k^{1*})^\top \mathbb{E}[y_{2,k+1} \omega_k \mid \mathcal{F}_{k-1}] \\
						&+ P_k x_k^* + M_k^\top v_k^* + p_k,
					\end{split}
				\end{equation}
				and $v^{*}_{k}=\Pi^{2*}_{k}x^{*}_{k}+\Sigma^{2*}_{k}$ in both \eqref{eq:3.7} and \eqref{eq:3.8}.
			\end{enumerate}
		\end{thm}

		\begin{proof}
			We establish necessity and sufficiency separately.
			
			\textbf{Necessity ($\Rightarrow$):}
			Assume quadruple strategy $\left(\Pi^{1*}, \Sigma^{1*}; \Pi^{2*}, \Sigma^{2*}\right)$ is a closed-loop Nash equilibrium. By Proposition \ref{prop:nash_equilibrium} (iii), for any $\mathbf{u} \in \mathcal{U}$,
			\begin{equation*}
				\mathcal{J}_1(0,\xi;\mathbf{u}^*,\mathbf{v}^*) \leq \mathcal{J}_1(0,\xi;\mathbf{u},\mathbf{v}^*)
			\end{equation*}
			where $\mathbf{v}^* = \Pi^{2*}\mathbf{x}^* + \Sigma^{2*}\in\mathcal{V}$ is given Player2's control input, $\mathbf{x}^*$ is generated by $(\mathbf{u}^*,\mathbf{v}^*)$, and $\mathbf{x}_u=\{\mathbf{x}_{u,k}\}_{k \in \mathcal{T}}$ is generated by $(\mathbf{u},\mathbf{v}^*)$.
			
			This implies $\mathbf{u}^*$ minimizes:
			\begin{equation}\label{eq:player1_problem}\begin{split}
					\tilde{\mathcal{J}}_1(0, \xi; \mathbf{u}, \mathbf{v}) =& \mathbb{E}\bigg[\frac{1}{2}\langle G_N x_{u,N},x_{u,N}\rangle + \langle g_N,x_{u,N}\rangle \\&+ \sum_{k=0}^{N-1} \ell_1(k,x_{u,k},u_k)\bigg],
				\end{split}
			\end{equation}
			subject to:
			\begin{equation}\nonumber
				\left\{\begin{aligned}
					x_{u,k+1} &= \tilde{A}_k x_{u,k} + B_k u_k + \tilde{b}_k + (\tilde{D}_k x_{u,k} + E_k u_k + \tilde{\sigma}_k)\omega_k,\\
					x_{u,0} &=\xi,
				\end{aligned}\right.
			\end{equation}
			where
			\begin{align*}
				\ell_1(k,x_{u,k},u_k) =& \frac{1}{2}\langle Q_k x_{u,k},x_{u,k}\rangle + \langle L_k x_{u,k},u_k\rangle
				\\&+ \frac{1}{2}\langle R_k u_k,u_k\rangle + \langle q_k,x_{u,k}\rangle + \langle \rho_k,u_k\rangle, \\ \nonumber
				\tilde{A}_k &= A_k + C_k \Pi_k^{2*}, \quad \tilde{b}_k = C_k \Sigma_k^{2*} + b_k, \\
				\tilde{D}_k &= D_k + F_k \Pi_k^{2*}, \quad \tilde{\sigma}_k = F_k \Sigma_k^{2*} + \sigma_k.
			\end{align*}
			By Theorem 3 in \cite{Wuyiwei2}, the first-order necessary conditions are exactly \eqref{eq:3.5} and \eqref{eq:3.6}.
			
			\textbf{Sufficiency ($\Leftarrow$):}
			Assume \eqref{eq:3.5}-\eqref{eq:3.8} hold. Under Assumption \ref{ass:2.3}:
			\begin{enumerate}
				\item [(i)] For Player 1: With fixed $\mathbf{v}^*\in\mathcal{V}$, \eqref{eq:3.5} and \eqref{eq:3.6} imply $\mathbf{u}^*$ is optimal for \eqref{eq:player1_problem}
				by Theorem 3 in \cite{Wuyiwei2}, so:
				\begin{equation*}
					\mathcal{J}_1(0,\xi;\mathbf{u}^*,\mathbf{v}^*) \leq \mathcal{J}_1(0,\xi;\mathbf{u},\mathbf{v}^*)
				\end{equation*}
				\item [(ii)] For Player 2: Similarly, \eqref{eq:3.7} and \eqref{eq:3.8} imply:
				\begin{equation*}
					\mathcal{J}_2(0,\xi;\mathbf{u}^*,\mathbf{v}^*) \leq \mathcal{J}_2(0,\xi;\mathbf{u}^*,\mathbf{v})
				\end{equation*}
			\end{enumerate}
			By Proposition \ref{prop:nash_equilibrium} (iii), the strategy is a closed-loop Nash equilibrium.
		\end{proof}

		\begin{rmk}
			~~~~~~\\
			The adjoint processes \(y_1 = \{y_{1,k}\}_{k\in \mathcal{T}}\) and \(y_2 = \{y_{2,k}\}_{k\in \mathcal{T}}\) in  Theorem \ref{thm:3.1A} are correspond to the feedback trajectories \((\mathbf{x}^*, \mathbf{u}^*, \mathbf{v}^*) = \left(\mathbf{x}^*, \Pi^{1*}\mathbf{x}^* + \Sigma^{1*}, \Pi^{2*}\mathbf{x}^* + \Sigma^{2*} \right)\).
		\end{rmk}
		Hamiltonian system characterizing Player 1's optimal control problem under another Player2's feedback information structure is
		\small{
			\begin{equation}\label{eq:2.22}
				\left\{\begin{array}{lll}
					\begin{split}
						x_{k+1}=&\big(A_{k}+C_{k}\Pi^{2}_{k}\big)x_{k}+B_{k}u_{k}+C_{k}\Sigma^{2}_{k}+b_{k}\\&+\big[\big(D_{k}+F_{k}\Pi^{2}_{k}\big)x_{k}+E_{k}u_{k}+F_{k}\Sigma^{2}_{k}+\sigma_{k}\big]\omega_{k},
						\\
						y_{1, k} =& \big(A_{k}+C_{k}\Pi^{2}_{k}\big)^{\top}\mathbb{E}[y_{1, k+1} \mid \mathcal{F}_{k-1}] + \big(D_{k}+F_{k}\Pi^{2}_{k}\big)^{\top}\\&\mathbb{E}[y_{1, k+1}\omega_{k} \mid \mathcal{F}_{k-1}] + Q_{k}x_{k} + L_{k}^{\top}u_{k}+q_{k},
						\\
						x_{0}=&\xi \in \mathbb R^n,~~~ y_{1, N}=G_{N}x_{N}+g_{N},
						\\&
						\frac{\partial \mathcal{H}_1(k, x_k, u_k, \Pi^{2}_{k}x_{k}+\Sigma^{2}_{k}, y_{1,k+1})}{\partial u_k}  \\
						=	& B_k^\top \mathbb{E}[y_{1,k+1}|\mathcal{F}_{k-1}] + E_k^\top \mathbb{E}[y_{1,k+1}\omega_k|\mathcal{F}_{k-1}] \\
						&+ L_k x_k + R_k u_k + \rho_k = 0,
					\end{split}
				\end{array}
				\right.
		\end{equation}}
		\normalsize{meanwhile,  Hamiltonian system for Player2 is}
		\small{
			\begin{equation}\label{eq:2.24}
				\left\{\begin{array}{lll}
					\begin{split}
						x_{k+1}=&\big(A_{k}+B_{k}\Pi^{1}_{k}\big)x_{k}+C_{k}v_{k}+B_{k}\Sigma^{1}_{k}+b_{k}\\&+\big[\big(D_{k}+E_{k}\Pi^{1}_{k}\big)x_{k}+F_{k}v_{k}+E_{k}\Sigma^{1}_{k}+\sigma_{k}\big]\omega_{k},
						\\
						y_{2, k} =& \big(A_{k}+B_{k}\Pi^{1}_{k}\big)^{\top}\mathbb{E}[y_{2, k+1} \mid \mathcal{F}_{k-1}] + \big(D_{k}+E_{k}\Pi^{1}_{k}\big)^{\top}\\&\mathbb{E}[y_{2, k+1}\omega_{k} \mid \mathcal{F}_{k-1}] + P_{k}x_{k} + M_{k}^{\top}v_{k}+p_{k},
						\\
						x_{0}=&\xi \in \mathbb R^n,~~~ y_{2, N} = H_{N}x_{N} + h_{N},
						\\
						&\frac{\partial \mathcal{H}_2(k, x_k, \Pi^{1}_{k}x_{k}+\Sigma^{1}_{k}, v_k, y_{2,k+1})}{\partial v_k}\\=&	C_{k}^{\top}\mathbb{E}[y_{2, k+1} \mid \mathcal{F}_{k-1}] + F_{k}^{\top}\mathbb{E}[y_{2, k+1}\omega_{k} \mid \mathcal{F}_{k-1}] \\&+ M_{k}x_{k}+ S_{k}v_{k}+\theta_{k} = 0,~~~k\in \mathcal{T}.
					\end{split}
				\end{array}
				\right.
		\end{equation}}
		\normalsize
		The Hamiltonian systems \eqref{eq:2.22} and \eqref{eq:2.24} represent the optimization conditions for each player, given the other player's strategy pair, and consist of two sets of fully coupled FBS$\Delta$Es.
		
		We postulate the following ansatzes for the relationship between the state process $\mathbf{x}$ and its adjoint processes $y_1 = \{y_{1,k}\}_{k\in\mathcal{T}}$ and $y_2 = \{y_{2,k}\}_{k\in\mathcal{T}}$:
		
		\begin{equation}\label{eq:2.25}
			\left\{\begin{array}{lll}
				\begin{split}
					y_{1,k} &= T^{1}_{k}x_{k} + \phi^{1}_{k}, \\
					y_{2,k} &= T^{2}_{k}x_{k} + \phi^{2}_{k},
				\end{split}
			\end{array}
			\right.
		\end{equation}
		with the terminal conditions given by:
		\begin{equation}\label{eq:3.13}
			\left\{\begin{array}{lll}
				\begin{split}
					T_N^1 &= G_N, & \phi_N^1 &= g_N, \\
					T_N^2 &= H_N, & \phi_N^2 &= h_N.
				\end{split}
			\end{array}
			\right.
		\end{equation}
		
		Futhermore, the above symmetric matrices $T^1$ and $T^2$ satisfy the uniform positive definiteness:
		$$
		T^1_k \succ 0 \quad \text{a.s.}, \quad \exists \delta>0 : T^1_k \succeq \delta I \quad \forall k \in \overline{\mathcal{T}}.
		$$
		$$
		T^2_k \succ 0 \quad \text{a.s.}, \quad \exists \delta>0 : T^2_k \succeq \delta I \quad \forall k \in \overline{\mathcal{T}}.
		$$
		
		Specifically:
		\begin{itemize}
			\item $T^{1} := \{T^{1}_{k}\}_{k=0}^N \in L_{\mathbb{F}}^{\infty}(\overline{\mathcal{T}}; \mathbb{S}_{+}^{n})$ with $T^{1}_{N} = G_{N}$.
			\item $T^{2} := \{T^{2}_{k}\}_{k=0}^N \in L_{\mathbb{F}}^{\infty}(\overline{\mathcal{T}}; \mathbb{S}_{+}^{n})$ with $T^{2}_{N} = H_{N}$.
			\item $\phi^{1} := \{\phi^{1}_{k}\}_{k=0}^N \in L_{\mathbb{F}}^2(\overline{\mathcal{T}}; \mathbb{R}^{n})$. satisfying the BS$\Delta$E:
			\begin{equation}\label{eq:2.26}
				\left\{\begin{array}{ll}
					\phi^{1}_{k} = g^{1}_{k+1}\mathbb{E}[\phi^{1}_{k+1} \mid \mathcal{F}_{k-1}] + h^{1}_{k+1}\mathbb{E}[\phi^{1}_{k+1}\omega_{k} \mid \mathcal{F}_{k-1}] + f^{1}_{k+1}, \\
					\phi^{1}_{N} = g_{N}.
				\end{array}\right.
			\end{equation}
			\item $\phi^{2} := \{\phi^{2}_{k}\}_{k=0}^N \in L_{\mathbb{F}}^2(\overline{\mathcal{T}}; \mathbb{R}^{n})$ satisfying the BS$\Delta$E:
			\begin{equation}\label{eq:2.27}
				\left\{\begin{array}{ll}
					\phi^{2}_{k} = g^{2}_{k+1}\mathbb{E}[\phi^{2}_{k+1} \mid \mathcal{F}_{k-1}] + h^{2}_{k+1}\mathbb{E}[\phi^{2}_{k+1}\omega_{k} \mid \mathcal{F}_{k-1}] + f^{2}_{k+1}, \\
					\phi^{2}_{N} = h_{N}.
				\end{array}\right.
			\end{equation}
		\end{itemize}
		
		Here, for i =1,2, the coefficients are $\mathcal{F}_{k-1}$-measurable processes:
		\[
		f^{i}_{k+1}  \in L_{\mathbb{F}}^2(\mathcal{T}; \mathbb{R}^n) , \quad
		g^{i}_{k+1}  \in L_{\mathbb{F}}^{\infty}(\mathcal{T}; \mathbb{R}^{n\times n}) ,\quad
		h^{i}_{k+1} \in L_{\mathbb{F}}^{\infty}(\mathcal{T}; \mathbb{R}^{n\times n}).
		\]
		
		Following the decoupling technique in Section 4.2 of \cite{Wuyiwei}, the stochastic Riccati equations are derived. The detailed similar derivations of CCREs are omitted here. Using the notations in Appendix A, for any $k\in\mathcal{T}$, the CCREs associated with DSDG-CL are formulated as follows:

		\begin{equation}\label{eq:2.28}
			\left\{\begin{array}{lll}
				\begin{split}
					T^{1}_{k}=&\Delta\left(T^{1}_{k+1},\Pi^{2}_{k}\right)
					-\Lambda\left(T^{1}_{k+1},\Pi^{2}_{k}\right)^{\top}\\&\times\Upsilon^{-1}(T^{1}_{k+1})\Lambda\left(T^{1}_{k+1},\Pi^{2}_{k}\right),
					\\
					\Pi^{1}_{k}=&-\Upsilon^{-1}(T^{1}_{k+1})\Lambda\left(T^{1}_{k+1},\Pi^{2}_{k}\right),
					\\
					T^{1}_{N}=&G_{N},
				\end{split}
			\end{array}
			\right.
		\end{equation}
		and
		\begin{equation}\label{eq:2.29}
			\left\{\begin{array}{lll}
				\begin{split}
					T^{2}_{k}
					=&\Delta\left(T^{2}_{k+1},\Pi^{1}_{k}\right)-	\Lambda\left(T^{2}_{k+1},\Pi^{1}_{k}\right)^{\top}\\&\times\Upsilon^{-1}(T^{2}_{k+1})	\Lambda\left(T^{2}_{k+1},\Pi^{1}_{k}\right),
					\\
					\Pi^{2}_{k}=&-\Upsilon^{-1}(T^{2}_{k+1})\Lambda\left(T^{2}_{k+1},\Pi^{1}_{k}\right),
					\\
					T^{2}_{N}=&H_{N}.
				\end{split}
			\end{array}
			\right.
		\end{equation}
		
		Furthermore, for any $k\in\mathcal{T}$, $\phi^{1}$ and $\phi^{2}$ satisfy the following BS$\Delta$Es respectively:
		\begin{equation}\label{eq:3.14}
			\left\{\begin{array}{lll}
				\begin{split}
					\phi^{1}_{k}=&\Theta\left(T^{1}_{k+1},\Pi^{2}_{k},\phi^{1}_{k+1}\right)-\Lambda^{\top}\left(T^{1}_{k+1},\Pi^{2}_{k}\right)	\\&\times\Upsilon^{-1}\left(T^{1}_{k+1}\right)\Phi\left(T^{1}_{k+1},\Pi^{2}_{k},\phi^{1}_{k+1}\right)
					,
					\\
					\Sigma^{1}_{k}=&-\Upsilon^{-1}\left(T^{1}_{k+1}\right)\Phi\left(T^{1}_{k+1},\Pi^{2}_{k},\phi^{1}_{k+1}\right),
					\\
					\phi^{1}_{N}=&g_{N},
				\end{split}
			\end{array}
			\right.
		\end{equation}
		and
		\begin{equation}\label{eq:3.15}
			\left\{\begin{array}{lll}
				\begin{split}		\phi^{2}_{k}=&\Theta\left(T^{2}_{k+1},\Pi^{1}_{k},\phi^{2}_{k+1}\right)-\Lambda^{\top}\left(T^{2}_{k+1},\Pi^{1}_{k}\right)\\&\times\Upsilon^{-1}\left(T^{2}_{k+1}\right)\Phi\left(T^{2}_{k+1},\Pi^{1}_{k},\phi^{2}_{k+1}\right)
					,
					\\
					\Sigma^{2}_{k}=&-\Upsilon^{-1}\left(T^{2}_{k+1}\right)\Phi\left(T^{2}_{k+1},\Pi^{1}_{k},\phi^{2}_{k+1}\right),
					\\
					\phi^{2}_{N}=&h_{N}.
				\end{split}
			\end{array}
			\right.
		\end{equation}.
		\begin{defn}\label{def:2.4}
			A set of solutions $(T^{1}, T^{2}) \in  L_{\mathbb{F}}^{\infty}(\overline{\mathcal{T}}; \mathbb{S}_{+}^{n})\times  L_{\mathbb{F}}^{\infty}
			(\overline{\mathcal{T}}; \mathbb{S}_{+}^{n}) $ of CCREs \eqref{eq:2.28} and \eqref{eq:2.29} is said to be regular if
			\begin{equation}\label{eq:3.16}
				\begin{array}{l}
					\left\{\begin{array}{l}
						\Upsilon(T^{1}_{k+1}),\Upsilon(T^{2}_{k+1}) \succeq 0,
						\\\\
						\mathcal{R}\big[\Lambda\left(T^{1}_{k+1},\Pi^{2}_{k}\right) \big]\subseteq \mathcal{R}\big[\Upsilon(T^{1}_{k+1})\big],
						\\\\
						\mathcal{R}\big[\Lambda\left(T^{2}_{k+1},\Pi^{1}_{k}\right)\big] \subseteq \mathcal{R}\big[\Upsilon(T^{2}_{k+1})\big],
						\\\\
						\Upsilon^{-1}(T^{1}_{k+1})\Lambda\left(T^{1}_{k+1},\Pi^{2}_{k}\right) \in \mathcal{Q}^{1}(\mathcal{T}),
						\\\\
						\Upsilon^{-1}(T^{2}_{k+1})\Lambda\left(T^{1}_{k+1},\Pi^{2}_{k}\right) \in \mathcal{Q}^{2}(\mathcal{T}).
					\end{array}\right.
				\end{array}
			\end{equation}
			The CCREs \eqref{eq:2.28}  and \eqref{eq:2.29} are said to be regularly solvable on $\mathcal{T}$ if they admit a set of regular solutions $(T^{1},T^{2})$.
		\end{defn}
		\begin{defn}\label{def:2.5}
		A set of solutions $(T^{1}, T^{2}) \in  L_{\mathbb{F}}^{\infty}(\overline{\mathcal{T}}; \mathbb{S}_{+}^{n})\times  L_{\mathbb{F}}^{\infty}
		(\overline{\mathcal{T}}; \mathbb{S}_{+}^{n}) $ of CCREs \eqref{eq:2.28} and \eqref{eq:2.29} is said to be strongly regular if there exists a constant $\delta > 0$ such that
		\begin{equation}\label{eq:3.21}
					\Upsilon(T^{1}_{k+1}),\Upsilon(T^{2}_{k+1}) \succ \delta I.
		\end{equation}
		\eqref{eq:2.28}  and \eqref{eq:2.29} are said to be strongly regularly solvable on $\mathcal{T}$ if they admit a set of strongly regular solutions $(T^{1},T^{2})$.
	\end{defn}
		Now, we present one of the main results on closed-loop Nash equilibrium for DSDG-CL: closed-loop state feedback representation.
		\begin{thm}\label{thm:4.3}
			\normalsize{Let Assumptions \ref{ass:2.1}-\ref{ass:2.3} hold.~Suppose two sets of CCREs. \eqref{eq:2.28} and
				\eqref{eq:2.29} admit
				a set of strongly regular solutions $(T^{1*}, T^{2*}, \Pi^{1*}, \Pi^{2*}) \in  L_{\mathbb{F}}^{\infty}(\mathcal{T}; \mathbb{S}_{+}^{n})\times  L_{\mathbb{F}}^{\infty}
				(\mathcal{T}; \mathbb{S}_{+}^{n})\times 	\mathcal{Q}^1(\mathcal{T})\times 	\mathcal{Q}^2(\mathcal{T}) $ and two sets of cross-coupled  BS$\Delta$Es \eqref{eq:3.14} and \eqref{eq:3.15} admit
				a set of solutions $(\phi^{1*}, \phi^{2*},\Sigma^{1*},\Sigma^{2*})\in  L_{\mathbb{F}}^2(\mathcal{T}; \mathbb{R}^{n})\times L_{\mathbb{F}}^2(\mathcal{T}; \mathbb{R}^{n})\times \mathcal{U}\times\mathcal{V}$.~Then Problem \ref{pro:1.2} is closed-loop solvable which closed-loop Nash equilibrium is $(\Pi^{1\ast}, \Sigma^{1\ast};\Pi^{2\ast}, \Sigma^{2\ast})\in  \mathscr{M}(\mathcal{T})$.And the corresponding closed-loop Nash equilibrium strategy pair $(\mathbf{u}^*,\mathbf{v}^*)$ has the following state feedback representation:}
			\small{	\begin{equation} \label{eq:4.1}
					\begin{split}
						&(\mathbf{u}^{\ast}, \mathbf{v}^{\ast})^{\top}\\=&(\Pi^{1\ast}\mathbf{x}^\ast+ \Sigma^{1\ast},\Pi^{2\ast}\mathbf{x}^\ast+ \Sigma^{2\ast})^{\top}\\=&-\left(\begin{array}{l}
							\Upsilon^{-1}(T^{1*}_{k+1})\left[\Lambda(T^{1*}_{k+1}, \Pi^{2*}_{k}) x^\ast_{k} + \Phi(T^{1*}_{k+1}, \Pi^{2*}_{k}, \phi^{1*}_{k+1})\right],
							\\\\
							\Upsilon^{-1}\left(T^{2*}_{k+1}\right)\left[\Lambda\left(T^{2*}_{k+1},\Pi^{1*}_{k}\right)x^\ast_{k}+\Phi\left(T^{2*}_{k+1},\Pi^{1*}_{k},\phi^{2*}_{k+1}\right)\right]
						\end{array}\right),
					\end{split}
			\end{equation}}
			\normalsize{where $\mathbf{x}^\ast \in L_{\mathbb{F}}^2(\overline{\mathcal{T}}; \mathbb{R}^n)$ is the solution of the following S$\Delta$E:}
			\begin{equation}\label{eq:4.2}
				\left\{
				\begin{array}{ll}
					\begin{split}
						x^{\ast}_{k+1} =& \big(A_{k}+B_{k} \Pi^{1\ast}_{k}+C_{k} \Pi^{2\ast}_{k}\big) x^{\ast}_{k} +B_{k}\Sigma^{1\ast}_{k} \\&+ C_{k} \Sigma^{2\ast}_{k} + b_{k} + \big[\big(D_{k}+ E_{k} \Pi^{1\ast}_{k}\\&+ F_{k} \Pi^{2\ast}_{k}\big)x^{\ast}_{k} +E_{k} \Sigma^{1\ast}_{k}  + F_{k} \Sigma^{2\ast}_{k}+\sigma_{k}\big] \omega_{k}, \\
						x_{0}^{\ast} =& \xi \in \mathbb{R}^{n}, \quad k \in \mathcal{T}.
					\end{split}
				\end{array}
				\right.
			\end{equation}
			
			\normalsize{Moreover, the corresponding Value function for Player1 is given by}
			\begin{equation}\label{eq:4.44}
				\begin{split}
					&V_1(0, \xi) \\= &\frac{1}{2} \mathbb{E} \sum_{k=0}^{N-1} \biggl\{ \bigg[2 \left\langle \phi^{1*}_{k+1}, (C_{k}\Sigma^{2*}_{k} + b_{k}) + (F_{k}\Sigma^{2*}_{k} + \sigma_{k}) \omega_{k} \right\rangle \\&- \left\langle \Phi(T^{1*}_{k+1}, \Pi^{2*}_{k}, \phi^{1*}_{k+1}), \Upsilon^{-1}(T^{1*}_{k+1}) \Phi(T^{1*}_{k+1}, \Pi^{2*}_{k}, \phi^{1*}_{k+1}) \right\rangle \\
					&+ \langle T^{1*}_{k+1} \big[(C_{k}\Sigma^{2*}_{k} + b_{k}) + (F_{k}\Sigma^{2*}_{k} + \sigma_{k}) \omega_{k}\big], (C_{k}\Sigma^{2*}_{k} + b_{k}) \\&+ (F_{k}\Sigma^{2*}_{k} + \sigma_{k}) \omega_{k} \rangle \bigg] + \left\langle \phi^{1*}_{0}, \xi \right\rangle + \left\langle T^{1*}_{0}\xi + \phi^{1*}_{0}, \xi \right\rangle \biggl\},
				\end{split}
			\end{equation}
			and for Player2:
			\begin{equation}\label{eq:4.45}
				\begin{split}
					&V_2(0,\xi)\\=&\frac{1}{2} \mathbb{E} \sum_{k=0}^{N-1} \biggl\{ \bigg[  2 \left\langle \phi^{2*}_{k+1}, (B_{k}\Sigma^{1*}_{k} + b_{k}) + (E_{k}\Sigma^{1*}_{k} + \sigma_{k}) \omega_{k} \right\rangle \\&- \left\langle \Phi(T^{2*}_{k+1}, \Pi^{1*}_{k}, \phi^{2*}_{k+1}), \Upsilon^{-1}(T^{2*}_{k+1}) \Phi(T^{2*}_{k+1}, \Pi^{1*}_{k}, \phi^{2*}_{k+1}) \right\rangle \\
					&+ \langle T^{2*}_{k+1} \left[(B_{k}\Sigma^{1*}_{k} + b_{k}) + (E_{k}\Sigma^{1*}_{k} + \sigma_{k}) \omega_{k}\right], (B_{k}\Sigma^{1*}_{k} + b_{k}) \\&+ (E_{k}\Sigma^{1*}_{k} + \sigma_{k}) \omega_{k} \rangle \bigg] + \left\langle \phi^{2*}_{0}, \xi \right\rangle + \left\langle T^{2*}_{0}\xi + \phi^{2*}_{0}, \xi \right\rangle \biggl\}.
				\end{split}
			\end{equation}
		\end{thm}
		\begin{proof}
			\subsubsection*{Step 1: Verification of state feedback representation of Closed-Loop Nash equilibrium}
			If the two sets of CCREs, \eqref{eq:2.28}--\eqref{eq:2.29}, admit a set of solutions $(T^{1*}, T^{2*}, \Pi^{1*}, \Pi^{2*}) $, and the two sets of cross-coupled BS$\Delta$Es, \eqref{eq:3.14}-\eqref{eq:3.15}, admit a set of solutions $(\phi^{1*}, \phi^{2*},\Sigma^{1*},\Sigma^{2*})$.~Using the notations \eqref{eq:2.30} in Appendix A,  next we prove the closed-loop Nash equilibrium $(\Pi^{1\ast}, \Sigma^{1\ast};\Pi^{2\ast}, \Sigma^{2\ast})$ having the state feedback representation \eqref{eq:4.1}, where $\mathbf{x}^\ast$ is the solution of the  S$\Delta$E \eqref{eq:4.2}.
			
			Define the following two Ansatzs for each Player:
\begin{equation}\label{eq:4.5}
y^\ast_{1,k}\equiv T^{1}_{k}x^\ast_{k}+\phi^{1*}_{k}
~~ k\in\mathcal{T},
\end{equation}

\begin{equation}\label{eq:4.6}
y^\ast_{2,k}\equiv T^{2}_{k}x^\ast_{k}+\phi^{2*}_{k}
~~ k\in\mathcal{T}.
\end{equation}
			
			From the two sets of CCREs.~\eqref{eq:2.28}-\eqref{eq:2.29}~and BS$\Delta$Es~\eqref{eq:3.14} - \eqref{eq:3.15}~presented earlier, we can derive $y^\ast_{1,k}$ and $y^\ast_{2,k}$ satisfies the following BS$\Delta$Es:
			\begin{equation}\label{eq:3.9}
				\left\{\begin{array}{lll}\begin{split}
						y^\ast_{1, k} =& \big(A_{k}+C_{k}\Pi^{2*}_{k}\big)^{\top}\mathbb{E}[y^\ast_{1, k+1} \mid \mathcal{F}_{k-1}] + \big(D_{k}+F_{k}\Pi^{2*}_{k}\big)^{\top}\\&\mathbb{E}[y^\ast_{1, k+1}\omega_{k} \mid \mathcal{F}_{k-1}] + Q_{k}x^\ast_{k} + L_{k}^{\top}u^\ast_{k}+q_{k},
						\\
						y^\ast_{1, N}=&G_{N}x^\ast_{N}+g_{N},~~~ k\in\mathcal{T},
					\end{split}
				\end{array}
				\right.
			\end{equation}
			and
			\begin{equation}\label{eq:3.10}
				\left\{\begin{array}{lll}
					\begin{split}
						y^\ast_{2, k} =& \big(A_{k}+B_{k}\Pi^{1*}_{k}\big)^{\top}\mathbb{E}[y^\ast_{2, k+1} \mid \mathcal{F}_{k-1}] + \big(D_{k}+E_{k}\Pi^{1*}_{k}\big)^{\top}\\&\mathbb{E}[y^\ast_{2, k+1}\omega_{k} \mid \mathcal{F}_{k-1}] + P_{k}x^\ast_{k} + M_{k}^{\top}v^\ast_{k}+p_{k},
						\\
						y^\ast_{2, N}=&H_{N}x^\ast_{N} + h_{N},~~~ k\in\mathcal{T}.
					\end{split}
				\end{array}
				\right.
			\end{equation}
			
			For Player~1, equation~\eqref{eq:3.9} gives the adjoint equation of the state equation~\eqref{eq:2.22} relative to the triple $\big(\Pi^{2*}, \Sigma^{2*}, \mathbf{x}^{*}\big)$, with a similar adjoint equation \eqref{eq:3.10} holding for Player~2, corresponding triple $\big(\Pi^{1*}, \Sigma^{1*}, \mathbf{x}^{*}\big)$.
			
			It is straightforward to verify that Player~1 satisfies the following stationarity condition:
			\begin{eqnarray}\label{eq:3.11}
				\begin{split}
					&B_{k}^{\top}\mathbb{E}[y^\ast_{1, k+1} \mid \mathcal{F}_{k-1}] + E_{k}^{\top}\mathbb{E}[y^\ast_{1, k+1}\omega_{k} \mid \mathcal{F}_{k-1}] \\&+ L_{k}x^\ast_{k} + R_{k}u^\ast_{k}+\rho_{k} = 0,~~~ k\in\mathcal{T}.
				\end{split}
			\end{eqnarray}
			Similarily for Player2,~we have
			\begin{eqnarray}\label{eq:4.53}
				\begin{split}
					&C_{k}^{\top}\mathbb{E}[y^\ast_{2, k+1} \mid \mathcal{F}_{k-1}] + F_{k}^{\top}\mathbb{E}[y^\ast_{2, k+1}\omega_{k} \mid \mathcal{F}_{k-1}] \\&+ M_{k}x^\ast_{k}+ S_{k}v^\ast_{k}+\theta_{k} = 0,~~~ k\in\mathcal{T}.
				\end{split}
			\end{eqnarray}
			Then, according to the sufficiency conditions in Theorem \ref{thm:3.1A}, it is clearly shown that $(\Pi^{1\ast}, \Sigma^{1\ast};\Pi^{2\ast}, \Sigma^{2\ast})$ is a closed-loop Nash equilibrium, and the corresponding closed-loop Nash equilibrium strategy pair $(\mathbf{u}^*,\mathbf{v}^*)$ satisfies the state feedback representation \eqref{eq:4.1}.
			\subsubsection*{Step 2: Value functions for each Player}
			To confirm that the Value functions match \eqref{eq:4.44} and \eqref{eq:4.45},~we give the following proof procedure for each Player.
			
			We start with Player 1. Given Player2's strategy $(\Pi^{2\ast},\Sigma^{2\ast})\in\mathscr{M}^{2}(\mathcal{T})$ and any $(\Pi^{1},\Sigma^{1})\in \mathscr{M}^{1}(\mathcal{T})$,~the total sum of the difference $\langle y_{1,k+1}^{\ast}, x_{k+1}^{(\Pi^{1},\Sigma^{1}, \Pi^{2\ast},\Sigma^{2\ast})}\rangle - \langle y_{1,k}^{\ast}, x_k^{(\Pi^{1},\Sigma^{1}, \Pi^{2\ast},\Sigma^{2\ast})}\rangle$  for each step $k$  from  $0$  through $N-1$ can be calculated that:
			\begin{equation}\label{eq:4.54}
				\begin{split}&
					\mathbb{E}\left\{\langle G_{N}x^{(\Pi^{1},\Sigma^{1}, \Pi^{2\ast},\Sigma^{2\ast})}_{N}+g_{N}, x_{N}^{(\Pi^{1},\Sigma^{1}, \Pi^{2\ast},\Sigma^{2\ast})} \rangle\right\} \\
					=& \mathbb{E}\left\{\langle y^\ast_{1, N}, x_{N}^{(\Pi^{1},\Sigma^{1}, \Pi^{2\ast},\Sigma^{2\ast})} \rangle\right\}  \\
					=& \mathbb{E}\sum_{k=0}^{N-1}\biggl\{\langle y^\ast_{1, k+1}, x_{k+1}^{(\Pi^{1},\Sigma^{1}, \Pi^{2\ast},\Sigma^{2\ast})} \rangle-\langle y^\ast_{1, k}, x_{k}^{(\Pi^{1},\Sigma^{1}, \Pi^{2\ast},\Sigma^{2\ast})} \rangle\\&+\langle y^\ast_{1, 0}, x_{0}^{(\Pi^{1},\Sigma^{1}, \Pi^{2\ast},\Sigma^{2\ast})} \rangle\biggr\}
					\\=& \mathbb{E}\sum_{k=0}^{N-1}\biggl\{ \bigg[\big\langle \big(A_{k}+C_{k}\Pi^{2*}_{k}\big)^{\top} y^\ast_{1, k+1} + \big(D_{k}+F_{k}\Pi^{2*}_{k}\big)^{\top}  y^\ast_{1, k+1} \omega_{k} \\& - y^\ast_{1, k}, x^{(\Pi^{1},\Sigma^{1}, \Pi^{2\ast},\Sigma^{2\ast})}_{k}\big\rangle + \big\langle B_{k}^{\top} y^\ast_{1, k+1} + E_{k}^{\top} y^\ast_{1, k+1} \omega_{k}, u^*_{k} \big\rangle \bigg]\\&+\big\langle y^\ast_{1, k+1},\big(C_{k}\Sigma^{2*}_{k}+b_{k}\big)+\big(F_{k}\Sigma^{2*}_{k}+\sigma_{k}\big)\omega_{k}\big\rangle\\&+\big\langle y^\ast_{1, 0},x^{(\Pi^{1},\Sigma^{1}, \Pi^{2\ast},\Sigma^{2\ast})}_{0}\big\rangle\biggr\}
					\\=& \mathbb{E}\sum_{k=0}^{N-1} \biggl\{\big[-\big\langle Q_{k}x^{(\Pi^{1},\Sigma^{1}, \Pi^{2\ast},\Sigma^{2\ast})}_{k} + L_{k}^{\top}u^\ast_{k} + q_{k}, x^{(\Pi^{1},\Sigma^{1}, \Pi^{2\ast},\Sigma^{2\ast})}_{k} \big\rangle\\& + \big\langle B_{k}^{\top} y^\ast_{1, k+1} + E_{k}^{\top} y^\ast_{1, k+1} \omega_{k}, u^*_{k} \big\rangle \big]+\big\langle T^{1*}_{k+1}x^{(\Pi^{1},\Sigma^{1}, \Pi^{2\ast},\Sigma^{2\ast})}_{k+1}\\&+\phi^{1*}_{k+1},\big(C_{k}\Sigma^{2*}_{k}+b_{k}\big)+\big(F_{k}\Sigma^{2*}_{k}+\sigma_{k}\big)\omega_{k}\big\rangle+\big\langle T^{1*}_{0}\xi+\phi^{1*}_{0},\xi\big\rangle\biggr\},
				\end{split}
			\end{equation}
			where the final equality is established based on the application of the adjoint equation of the \eqref{eq:4.5} and \eqref{eq:3.9} .
			
			By substituting   \eqref{eq:4.54} into the corresponding cost functional  $\mathcal{J}_{1}(0,\xi;\Pi^{1},\Sigma^{1},\Pi^{2\ast},\Sigma^{2\ast})$ defined by \eqref{eq:1.4}, we get
			\begin{equation}\nonumber
				\begin{split}
					\mathcal{J}_1&(0, \xi, \Pi^{1},\Sigma^{1}, \Pi^{2\ast},\Sigma^{2\ast})\\=&\frac{1}{2} \mathbb{E}\bigg\{\big\langle G_{N} x^{(\Pi^{1},\Sigma^{1}, \Pi^{2\ast},\Sigma^{2\ast})}_{N}, x^{(\Pi^{1},\Sigma^{1}, \Pi^{2\ast},\Sigma^{2\ast})}_{N}\big\rangle \\&+ 2\big\langle g_{N}, x^{(\Pi^{1},\Sigma^{1}, \Pi^{2\ast},\Sigma^{2\ast})}_{N}\big\rangle \\
					&+ \sum_{k=0}^{N-1} \bigg[ \big\langle Q_{k}x^{(\Pi^{1},\Sigma^{1}, \Pi^{2\ast},\Sigma^{2\ast})}_{k}, x^{(\Pi^{1},\Sigma^{1}, \Pi^{2\ast},\Sigma^{2\ast})}_{k}\big\rangle\\&+ 2\big\langle L_{k}^{\top} u^\ast_{k}, x^{(\Pi^{1},\Sigma^{1}, \Pi^{2\ast},\Sigma^{2\ast})}_{k}\big\rangle + \langle R_{k}u^\ast_{k}, u^\ast_{k}\rangle\\&  + 2\big\langle q_{k}, x^{(\Pi^{1},\Sigma^{1}, \Pi^{2\ast},\Sigma^{2\ast})}_{k}\big\rangle + 2\langle \rho_{k}, u^\ast_{k}\rangle\bigg]\bigg\}
				\end{split}
			\end{equation}
			\begin{equation}\label{eq:4.55}
				\begin{split}=&\frac{1}{2} \mathbb{E}\sum_{k=0}^{N-1}\biggl\{ \bigg[\big\langle B_{k}^{\top}y^{(\Pi^{1},\Sigma^{1}, \Pi^{2\ast},\Sigma^{2\ast})}_{1, k+1}  + E_{k}^{\top}\\
					&\times y^{(\Pi^{1},\Sigma^{1}, \Pi^{2\ast},\Sigma^{2\ast})}_{1, k+1}\omega_{k}  + L_{k}x^{(\Pi^{1},\Sigma^{1}, \Pi^{2\ast},\Sigma^{2\ast})}_{k} \\
					&+ R_{k}u^\ast_{k}+\rho_{k},u^\ast_{k}\big\rangle+\big\langle q_{k},x^{(\Pi^{1},\Sigma^{1}, \Pi^{2\ast},\Sigma^{2\ast})}_{k}\big\rangle\\&+\langle \rho_{k},u^\ast_{k}\rangle+\big\langle g_{N} , x^{(\Pi^{1},\Sigma^{1}, \Pi^{2\ast},\Sigma^{2\ast})}_{N}\big\rangle+\big\langle T^{1*}_{k+1}\\&\times x^{(\Pi^{1},\Sigma^{1}, \Pi^{2\ast},\Sigma^{2\ast})}_{k+1}+\phi^{1*}_{k+1},\big(C_{k}\Sigma^{2*}_{k}+b_{k}\big)\\&+\big(F_{k}\Sigma^{2*}_{k}+\sigma_{k}\big)\omega_{k}\big\rangle\bigg]+\langle T^{1*}_{0}\xi+\phi^{1*}_{0},\xi\rangle\biggl\}
					\\=&\frac{1}{2} \mathbb{E}\sum_{k=0}^{N-1}\biggl\{ \bigg[\big\langle q_{k},x^{(\Pi^{1},\Sigma^{1}, \Pi^{2\ast},\Sigma^{2\ast})}_{k}\big\rangle+\langle \rho_{k},u^\ast_{k}\rangle\\
					&+\big\langle g_{N} , x^{(\Pi^{1},\Sigma^{1}, \Pi^{2\ast},\Sigma^{2\ast})}_{N}\big\rangle+\big\langle T^{1*}_{k+1}x^{(\Pi^{1},\Sigma^{1}, \Pi^{2\ast},\Sigma^{2\ast})}_{k+1}\\	&+\phi^{1*}_{k+1},\big(C_{k}\Sigma^{2*}_{k}+b_{k}\big)+\big(F_{k}\Sigma^{2*}_{k}+\sigma_{k}\big)\omega_{k}\big\rangle\bigg]\\&+\langle T^{1*}_{0}\xi+\phi^{1*}_{0},\xi\rangle\biggl\},
				\end{split}
			\end{equation}
			where the stationarity condition \eqref{eq:3.11} is employed to prove the final equation.~Next, we furthermore arrive at the following equation of $\mathbb E\big[\big\langle g_{N} , x^{(\Pi^{1},\Sigma^{1}, \Pi^{2\ast},\Sigma^{2\ast})}_{N}\big\rangle\big]$ by aggregating $\big\langle \phi^{1*}_{k+1}, x^{(\Pi^{1},\Sigma^{1}, \Pi^{2\ast},\Sigma^{2\ast})}_{k+1}\big\rangle - \big\langle \phi^{1*}_{k}, x^{(\Pi^{1},\Sigma^{1}, \Pi^{2\ast},\Sigma^{2\ast})}_{k}\big\rangle$ from $k=0$ to $k=N$:
			\begin{equation}\label{eq:4.56}
				\begin{split}
					&\mathbb E[\langle g_{N} , x^{(\Pi^{1},\Sigma^{1}, \Pi^{2\ast},\Sigma^{2\ast})}_{N}\rangle]\\=&\mathbb{E}\sum_{k=0}^{N-1}\biggl\{\bigg[\big\langle \phi^{1*}_{k+1}, x^{(\Pi^{1},\Sigma^{1}, \Pi^{2\ast},\Sigma^{2\ast})}_{k+1}\big\rangle \\&- \big\langle \phi^{1*}_{k}, x^{(\Pi^{1},\Sigma^{1}, \Pi^{2\ast},\Sigma^{2\ast})}_{k}\big\rangle\bigg]\biggr\}+\big\langle \phi^{1*}_{0}, x^{(\Pi^{1},\Sigma^{1}, \Pi^{2\ast},\Sigma^{2\ast})}_{0}\big\rangle.
				\end{split}
			\end{equation}
			Putting \eqref{eq:4.56} into \eqref{eq:4.55},  we have
			\begin{equation}\nonumber
				\begin{split}
					\mathcal{J}&_1(0, \xi, \Pi^{1},\Sigma^{1}, \Pi^{2\ast},\Sigma^{2\ast})\\=&\frac{1}{2} \mathbb{E}\sum_{k=0}^{N-1}\biggl\{ \bigg[\langle q_{k},x^{(\Pi^{1},\Sigma^{1}, \Pi^{2\ast},\Sigma^{2\ast})}_{k}\rangle+\langle \rho_{k},u^\ast_{k}\rangle\\&+\langle \phi^{1*}_{k+1}, x^{(\Pi^{1},\Sigma^{1}, \Pi^{2\ast},\Sigma^{2\ast})}_{k+1}\rangle+\big\langle T^{1*}_{k+1}x^{(\Pi^{1},\Sigma^{1}, \Pi^{2\ast},\Sigma^{2\ast})}_{k+1}\\&+\phi^{1*}_{k+1},\big(C_{k}\Sigma^{2*}_{k}+b_{k}\big)+\big(F_{k}\Sigma^{2*}_{k}+\sigma_{k}\big)\omega_{k}\big\rangle\\&- \langle \phi^{1*}_{k}, x^{(\Pi^{1},\Sigma^{1}, \Pi^{2\ast},\Sigma^{2\ast})}_{k}\rangle\bigg] +\langle \phi^{1*}_{0}, \xi\rangle+\langle T^{1*}_{0}\xi,\xi\rangle\biggl\}
					\\	=&\frac{1}{2} \mathbb{E}\sum_{k=0}^{N-1}\biggl\{ \bigg[\big\langle q_{k},x^{(\Pi^{1},\Sigma^{1}, \Pi^{2\ast},\Sigma^{2\ast})}_{k}\big\rangle-\bigg\langle \rho_{k},\Upsilon^{-1}\left(T^{1*}_{k+1}\right)\\&\times\bigg[\Lambda\left(T^{2*}_{k+1},\Pi^{1*}_{k}\right) x^\ast_{k}+\Phi\left(T^{1*}_{k+1},\Pi^{2*}_{k},\phi^{1*}_{k+1}\right)\bigg]\bigg\rangle\end{split}
			\end{equation}
			\begin{equation}\label{eq:4.16}
				\begin{split}&+\bigg\langle \phi^{1*}_{k+1}, \big(A_{k}+C_{k}\Pi^{2*}_{k}\big)x^\ast_{k}-B_{k}\biggl\{\Upsilon^{-1}\left(T^{1*}_{k+1}\right)\\&\times\bigg[\Lambda\left(T^{2*}_{k+1},\Pi^{1*}_{k}\right) x^\ast_{k}+\Phi\left(T^{1*}_{k+1},\Pi^{2*}_{k},\phi^{1*}_{k+1}\right)\bigg]\biggr\}
					\\&+C_{k}\Sigma^{2*}_{k}+b_{k}+\bigg[\big(D_{k}+F_{k}\Pi^{2*}_{k}\big)x^\ast_{k}-E_{k}\biggl\{\Upsilon^{-1}\left(T^{1*}_{k+1}\right)
					\\
					&\times\bigg[\Lambda\left(T^{2*}_{k+1},\Pi^{1*}_{k}\right)x^\ast_{k}+\Phi\left(T^{1*}_{k+1},\Pi^{2*}_{k},\phi^{1*}_{k+1}\right)\bigg]\biggr\}\\&+F_{k}\Sigma^{2*}_{k}+\sigma_{k}\bigg]\omega_{k}\bigg\rangle+\bigg\langle T^{1*}_{k+1}\bigg[\big(A_{k}+C_{k}\Pi^{2*}_{k}\big)x^\ast_{k}-B_{k}\\&\times\biggl\{\Upsilon^{-1}\left(T^{1*}_{k+1}\right)\bigg[\Lambda\left(T^{2*}_{k+1},\Pi^{1*}_{k}\right) x^\ast_{k}+\Phi\left(T^{1*}_{k+1},\Pi^{2*}_{k},\phi^{1*}_{k+1}\right)\bigg]	\biggr\}\\&+C_{k}\Sigma^{2*}_{k}+b_{k}+\bigg[\big(D_{k}+F_{k}\Pi^{2*}_{k}\big) x^\ast_{k}-E_{k}\biggl\{\Upsilon^{-1}\left(T^{1*}_{k+1}\right)\\&\times\bigg[\Lambda\left(T^{2*}_{k+1},\Pi^{1*}_{k}\right)x^\ast_{k}+\Phi\left(T^{1*}_{k+1},\Pi^{2*}_{k},\phi^{1*}_{k+1}\right)\bigg]\biggr\}+F_{k}\Sigma^{2*}_{k}\\&+\sigma_{k}\bigg]\omega_{k}\bigg]+\phi^{1*}_{k+1},\big(C_{k}\Sigma^{2*}_{k}+b_{k}\big)+\big(F_{k}\Sigma^{2*}_{k}+\sigma_{k}\big)\omega_{k}\bigg\rangle\\&- \big\langle \phi^{1*}_{k}, x^{(\Pi^{1},\Sigma^{1}, \Pi^{2\ast},\Sigma^{2\ast})}_{k}\big\rangle\bigg]+\big\langle \phi^{1*}_{0}, \xi\big\rangle +\big\langle T^{1*}_{0}\xi,\xi\big\rangle\biggl\}
					\\=&\frac{1}{2} \mathbb{E}\sum_{k=0}^{N-1}\biggl\{\bigg[\big\langle \Theta\left(T^{1*}_{k+1},\Pi^{2*}_{k},\phi^{1*}_{k+1}\right)-\Lambda^{\top}\left(T^{1*}_{k+1},\Pi^{2*}_{k}\right)\big[\Upsilon^{-1}\left(T^{1*}_{k+1}\right)\\&\times\Phi\left(T^{1*}_{k+1},\Pi^{2*}_{k},\phi^{1*}_{k+1}\right)\big]
					-\phi^{1*}_{k},x^{(\Pi^{1},\Sigma^{1}, \Pi^{2\ast},\Sigma^{2\ast})}_{k}\big\rangle+2\big\langle\phi^{1*}_{k+1},\\&\big(C_{k}\Sigma^{2*}_{k}+b_{k}\big)+\big(F_{k}\Sigma^{2*}_{k}+\sigma_{k}\big)\omega_{k}\big\rangle-\big\langle\Phi\left(T^{1*}_{k+1},\Pi^{2*}_{k},\phi^{1*}_{k+1}\right),\\&\Upsilon^{-1}\left(T^{1*}_{k+1}\right)\Phi\left(T^{1*}_{k+1},\Pi^{2*}_{k},\phi^{1*}_{k+1}\right)\big\rangle+\big\langle T^{1*}_{k+1}\big[\big(C_{k}\Sigma^{2*}_{k}+b_{k}\big)\\&+\big(F_{k}\Sigma^{2*}_{k}+\sigma_{k}\big)\omega_{k}\big],\big(C_{k}\Sigma^{2*}_{k}+b_{k}\big)+\big(F_{k}\Sigma^{2*}_{k}+\sigma_{k}\big)\omega_{k}\big\rangle\bigg]\\&+\big\langle \phi^{1*}_{0}, \xi\big\rangle+\big\langle T^{1*}_{0}\xi,\xi\big\rangle\biggl\},
				\end{split}
			\end{equation}
			where the notations in Appendix A and the BS$\Delta$Es \eqref{eq:3.9} were employed to prove the final equality.
			
			Let Assumptions \ref{ass:2.1}-\ref{ass:2.2} hold, it follows immediately that
			\begin{equation}\label{eq:4.58}
				\begin{split}
					\mathcal{J}_1(0,\xi;\Pi^{1},\Sigma^{1},&\Pi^{2\ast},\Sigma^{2\ast})\geq\mathcal{J}_1(0,\xi;\Pi^{1\ast},\Sigma^{1\ast},\Pi^{2\ast},\Sigma^{2\ast}),\\\ &\qquad\qquad\qquad\forall (\Pi^{1},\Sigma^{1})\in \mathscr{M}^{1}(\mathcal{T}).
				\end{split}
			\end{equation}
			In \eqref{eq:4.58}, the equality holds if and only if $(\Pi^{1},\Sigma^{1})=(\Pi^{1\ast},\Sigma^{1\ast})$, that is, $(\Pi^{1\ast},\Sigma^{1\ast})$ is  a minimizer of $\mathcal{J}_1(0,\xi;\Pi^{1},\Sigma^{1},\Pi^{2\ast},\Sigma^{2\ast})$.~Thus,~we can obtain the Value function for Player1:
			\begin{equation}	\label{eq:4.59}
				\begin{split}
					V&_1(0, \xi) \\=&\mathcal{J}_1(0, \xi, \Pi^{1\ast}\mathbf{x}^{\ast}+\Sigma^{1\ast}, \Pi^{2\ast}\mathbf{x}^{\ast}+\Sigma^{2\ast})
					\\=&\frac{1}{2} \mathbb{E} \sum_{k=0}^{N-1} \biggl\{ \bigg[2 \left\langle \phi^{1*}_{k+1}, (C_{k}\Sigma^{2*}_{k} + b_{k}) + (F_{k}\Sigma^{2*}_{k} + \sigma_{k}) \omega_{k} \right\rangle \\&- \left\langle \Phi(T^{1*}_{k+1}, \Pi^{2*}_{k}, \phi^{1*}_{k+1}), \Upsilon^{-1}(T^{1*}_{k+1}) \Phi(T^{1*}_{k+1}, \Pi^{2*}_{k}, \phi^{1*}_{k+1}) \right\rangle \\
					&+ \langle T^{1*}_{k+1} \big[(C_{k}\Sigma^{2*}_{k} + b_{k}) + (F_{k}\Sigma^{2*}_{k} + \sigma_{k}) \omega_{k}\big], (C_{k}\Sigma^{2*}_{k} + b_{k}) \\&+ (F_{k}\Sigma^{2*}_{k} + \sigma_{k}) \omega_{k} \rangle \bigg] + \left\langle \phi^{1*}_{0}, \xi \right\rangle + \left\langle T^{1*}_{0}\xi + \phi^{1*}_{0}, \xi \right\rangle \biggl\}.
				\end{split}
			\end{equation}
			Utilizing analogous techniques applied to Player 1, we can derive the Value function expression for Player 2:
			\begin{equation}\label{eq:4.65}
				\begin{split}
					V&_2(0, \xi) \\=&\frac{1}{2} \mathbb{E} \sum_{k=0}^{N-1} \biggl\{ \bigg[  2 \left\langle \phi^{2*}_{k+1}, (B_{k}\Sigma^{1*}_{k} + b_{k}) + (E_{k}\Sigma^{1*}_{k} + \sigma_{k}) \omega_{k} \right\rangle \\&- \left\langle \Phi(T^{2*}_{k+1}, \Pi^{1*}_{k}, \phi^{2*}_{k+1}), \Upsilon^{-1}(T^{2*}_{k+1}) \Phi(T^{2*}_{k+1}, \Pi^{1*}_{k}, \phi^{2*}_{k+1}) \right\rangle \\
					&+ \langle T^{2*}_{k+1} \left[(B_{k}\Sigma^{1*}_{k} + b_{k}) + (E_{k}\Sigma^{1*}_{k} + \sigma_{k}) \omega_{k}\right], (B_{k}\Sigma^{1*}_{k} + b_{k}) \\&+ (E_{k}\Sigma^{1*}_{k} + \sigma_{k}) \omega_{k} \rangle \bigg] + \left\langle \phi^{2*}_{0}, \xi \right\rangle + \left\langle T^{2*}_{0}\xi + \phi^{2*}_{0}, \xi \right\rangle \biggl\}.
				\end{split}
			\end{equation}
			
			To sum up, the existence of the closed-loop Nash equilibrium  $(\Pi^{1\ast}, \Sigma^{1\ast};\Pi^{2\ast}, \Sigma^{2\ast})$ is an immediate result by combining \eqref{eq:4.59} and \eqref{eq:4.65}. Furthermore, the closed-loop Nash equilibrium is specified by \eqref{eq:4.1}. Problem \ref{pro:1.2} admits a linear state feedback representation via the solution to the two sets of  CCREs.\eqref{eq:2.28}-\eqref{eq:2.29}. The proof is complete.
		\end{proof}

\begin{rmk}\label{rem:riccati-existence}
\textbf{(Solvability of Riccati Equations and Open Challenges)}

Theorem \ref {thm:4.3} assumes the existence of a strongly regular solution $(T^{1*}, T^{2*}, \Pi^{1*}, \Pi^{2*})$ to the CCREs \eqref{eq:2.28} and \eqref{eq:2.29} . The justification for this assumption stems from the following observations:

\begin{enumerate}
    \item [(i)]\textbf{Degenerated Case (Optimal Control Problem):}
    When the problem degenerates to a single-player optimal control problem (a special case of non-cooperative games), explicit solutions to the Riccati equations can be constructed via backward induction under Assumption \ref{ass:2.3} . The noise independence and condition $R_k \succ \delta I$ jointly guarantee the boundedness and invertibility required for the iterative process.

    \item [(ii)]\textbf{Random Coefficient Challenge:}
    In non-zero-sum games with random coefficients, the coupling in CCREs (i.e., $\Pi_k^{1*}$ depending on $T_{k+1}^2$, and $\Pi_k^{2*}$ depending on $T_{k+1}^1$) renders solution existence highly nontrivial. Existing theory has only established existence proofs for special structures, such as diagonally dominant coupling terms or decoupled noise. The solvability of general-form stochastic CCREs remains an \textbf{open problem} requiring further investigation.
\end{enumerate}
\end{rmk}

		\section{CCRE's Closed-Loop Solvability}\label{sec:4}
		By the Remark 2.5.4 in \cite{sun4} and game theory in \cite{sun2020game}, the related assumptions ensure the strongly regular solvability of Riccati equations, implying open- or closed-loop solvability in stochastic LQ optimal control problems or stochastic LQ non-zero-sum differential/difference games, which can be regarded as two stochastic LQ optimal control problems (SLQs,~for short). Factually, when weight matrices $R_k$ or $S_k$ are non-positive-definite (singular), the Nash equilibrium may still exist under specific structural conditions (e.g., decoupled dynamics or collaborative cost couplings). In this section, we discuss the closed-loop solvability of DSDG-CL  when the Assumption \ref{ass:2.3} is invalid.
		
		To further develop our analysis, the next proposition provides an essential theoretical foundation for performing a coherent analysis.
		
		\begin{prop}\label{prop4.1}
			Let Assumptions \ref{ass:2.1}-\ref{ass:2.2} hold. If $\left(\Pi^{1*}, \Sigma^{1*}; \Pi^{2*}, \Sigma^{2*}\right) \in \mathscr{M}(\mathcal{T})$ constitutes a closed-loop Nash equilibrium for the DSDG-CL over $\mathcal{T}$ and let $\mathbb{X}(\cdot)$ be the solution to the $\mathbb{R}^{n\times n}$-valued S$\Delta$Es
				\begin{equation}\label{eq:4.1}
				\left\{
				\begin{aligned}
					\mathbb{X}^*_{k+1} &= \big(A_{k}+B_{k}\Pi^{1*}_{k}+C_{k}\Pi^{2*}_{k}\big)\mathbb{X}^*_{k} \\
					&\quad + \big[\big(D_{k}+E_{k}\Pi^{1*}_{k}+F_{k}\Pi^{2*}_{k}\big)\mathbb{X}^*_{k}\big]\omega_{k}, \\
					\mathbb{X}^*_{0} &= I, \quad k\in \mathcal{T},
				\end{aligned}
				\right.
			\end{equation}
			then the following two stationarity conditions for each Player hold simultaneously:
				\begin{enumerate}
				\item For Player 1, the adapted solution $\mathbb{Y}^*_{1}$ to the $\mathbb{R}^{n \times n}$-valued BS$\Delta$E
					\begin{equation}\label{eq:4.2}
					\left\{
					\begin{aligned}
						\mathbb{Y}^*_{1, k} &= \big(A_{k}+C_{k}\Pi^{2*}_{k}\big)^{\top}\mathbb{E}[\mathbb{Y}^*_{1, k+1} \mid \mathcal{F}_{k-1}] \\
						&\quad + \big(D_{k}+F_{k}\Pi^{2*}_{k}\big)^{\top}\mathbb{E}[\mathbb{Y}^*_{1, k+1}\omega_{k} \mid \mathcal{F}_{k-1}] \\
						&\quad + \big(Q_{k}+L_{k}^{\top}\Pi^{1*}_{k}\big)\mathbb{X}^*_{k} ,  \\\mathbb{Y}^*_{1, N}&=G_{N}\mathbb{X}^*_{N}, \quad k\in \mathcal{T},
					\end{aligned}
					\right.
				\end{equation}
			satisfies
				\begin{equation}\label{eq:4.3}
				\begin{split}
					&B_{k}^{\top}\mathbb{E}[\mathbb{Y}^*_{1, k+1} \mid \mathcal{F}_{k-1}] + E_{k}^{\top}\mathbb{E}[\mathbb{Y}^*_{1, k+1}\omega_{k} \mid \mathcal{F}_{k-1}] \\
					&\quad + \big(L_{k}+R_{k}\Pi^{1*}_{k}\big)\mathbb{X}^*_{k} = 0, \quad k\in \mathcal{T}.
				\end{split}
			\end{equation}
				\item For Player 2, the adapted solution $\mathbb{Y}^*_{1}$ to the $\mathbb{R}^{n \times n}$-valued BS$\Delta$E
				\begin{equation}\label{eq:4.4}
					\left\{
					\begin{aligned}
						\mathbb{Y}^*_{2, k} &= \big(A_{k}+B_{k}\Pi^{1*}_{k}\big)^{\top}\mathbb{E}[\mathbb{Y}^*_{2, k+1} \mid \mathcal{F}_{k-1}] \\
						&\quad + \big(D_{k}+E_{k}\Pi^{1*}_{k}\big)^{\top}\mathbb{E}[\mathbb{Y}^*_{2, k+1}\omega_{k} \mid \mathcal{F}_{k-1}] \\
						&\quad + \big(P_{k}+M_{k}^{\top}\Pi^{2*}_{k}\big)\mathbb{X}^*_{k}, \\
						 \mathbb{Y}^*_{2, N}&=H_{N}\mathbb{X}^*_{N}, \quad k\in \mathcal{T},
					\end{aligned}
					\right.
				\end{equation}
			satisfies
				\begin{equation}\label{eq:4.5}
				\begin{split}
					&C_{k}^{\top}\mathbb{E}[\mathbb{Y}^{\ast}_{2, k+1} \mid \mathcal{F}_{k-1}] + F_{k}^{\top}\mathbb{E}[\mathbb{Y}^{\ast}_{2, k+1}\omega_{k} \mid \mathcal{F}_{k-1}] \\
					&\quad + \big(M_{k}+S_k\Pi^{2*}_{k}\big)\mathbb{X}^{\ast}_{k} = 0, \quad k\in \mathcal{T}.
				\end{split}
			\end{equation}
				\end{enumerate}
		\end{prop}
		
		\begin{proof}
			Consider the state equation:
			\begin{equation}\label{eq:40.1}
				\begin{cases}
					\begin{split}
						x_{k+1} &= \big(A_k+B_k\Pi^{1*}_{k}+C_k\Pi^{2*}_{k}\big)x_k + B_k \Sigma^{1}_{k} + C_k \Sigma^{2*}_{k} \\
						&\quad + b_k + \big[\big(D_k+E_k\Pi^{1*}_{k}+F_k\Pi^{2*}_{k}\big)x_k + E_k \Sigma^{1}_{k} \\
						&\quad + F_k \Sigma^{2*}_{k} + \sigma_k\big]\omega_k \\
						x_0 &= \xi.
					\end{split}
				\end{cases}
			\end{equation}
			with the cost functionals for each Player defined as:
			\begin{equation}\label{eq:40.2}
				\begin{split}
					\tilde{\mathcal{J}}_1(0, \xi; \Sigma^{1}) \equiv& \mathcal{J}_1(0, \xi; \Pi^{1*}\mathbf{x}+\Sigma^{1}, \Pi^{2*}\mathbf{x}+\Sigma^{2*})
					\\=&\frac{1}{2} \mathbb{E}\bigg\{\langle G_{N} x_{N}, x_{N}\rangle + 2\langle g_{N}, x_{N}\rangle \\
					& + \sum_{k=0}^{N-1} \bigg[\langle \big[Q_{k}+ L_{k}^{\top}\Pi^{1*}_{k}+ \Pi^{1*\top}_{k}L_{k}\\
					&+\Pi^{1*\top}_{k}R_{k}\Pi^{1*}_{k}\big]x_{k}, x_{k}\rangle  +2\langle \big[ L_{k}^{\top}\Sigma^{1}_{k}\\
					&+\Pi^{1*\top}_{k}R_{k}\Sigma^{1}_{k}+\rho_{k}\Pi^{1*}_{k}+q_{k}\big], x_{k}\rangle  \\
					&+\langle R_{k}\Sigma^{1}_{k}, \Sigma^{1}_{k}\rangle+2\langle \rho_{k}, \Sigma^{1}_{k}\rangle\bigg]\bigg\},
				\end{split}
			\end{equation}
			\begin{equation}\label{eq:40.3}
			\begin{split}
				\tilde{\mathcal{J}}_2(0, \xi; \Sigma^{2}) \equiv& \mathcal{J}_2(0, \xi; \Pi^{1*}\mathbf{x}+\Sigma^{1*}, \Pi^{2*}\mathbf{x}+\Sigma^{2})
				\\=&\frac{1}{2} \mathbb{E}\bigg\{\langle H_{N} x_{N}, x_{N}\rangle + 2\langle h_{N}, x_{N}\rangle \\
				& +\sum_{k=0}^{N-1} \bigg[\langle \big[P_{k}+M^{\top}_{k}\Pi^{2*}_{k}+\Pi^{2*\top}_{k}M_{k}\\
				&+\Pi^{2*\top}_{k}S_{k}\Pi^{2*}_{k}\big]x_{k}, x_{k}\rangle  + 2\langle\big[ M_{k}^{\top}\Sigma^{2}_{k}\\
				&+\Pi^{2*\top}_{k}S_{k}\Sigma^{2}_{k}+\theta_{k}\Pi^{2*}_{k}+p_{k}\big], x_{k}\rangle \\
				& + \langle S_{k}\Sigma^{2}_{k}, \Sigma^{2}_{k}\rangle + 2\langle \theta_{k}, \Sigma^{2}_{k}\rangle\bigg]\bigg\}.
			\end{split}
		\end{equation}
			
			By the implication (i) $\Rightarrow$ (ii) of Proposition \ref{prop:nash_equilibrium}, for any initial state $\xi\in\mathbb{R}^{n}$, if $\left(\Pi^{1*}, \Sigma^{1*}; \Pi^{2*}, \Sigma^{2*}\right) \in \mathscr{M}(\mathcal{T})$ is a closed-loop Nash equilibrium for DSDG-CL over $\mathcal{T}$, then $(\Sigma^{1*},\Sigma^{2*})\in\mathcal{U}\times\mathcal{V}$ constitutes an open-loop Nash equilibrium for the problem with state equation \eqref{eq:40.1} and cost functionals \eqref{eq:40.2} and \eqref{eq:40.3}.
			
			According to Theorem 4.1 in \cite{Wuyiwei}, we obtain the necessary conditions:
			\begin{equation}\label{eq:40.4}
				\begin{split}
					&B_k^\top \mathbb{E}[y^*_{1,k+1} \mid \mathcal{F}_{k-1}] + E_k^\top \mathbb{E}[y^*_{1,k+1} \omega_k \mid \mathcal{F}_{k-1}] \\
					&\quad + (L_k+R_k\Pi^{1*}_{k}) x^{*}_k + R_k \Sigma^{1}_k + \rho_k = 0,
				\end{split}
			\end{equation}
			
			\begin{equation}\label{eq:40.5}
				\begin{split}
					&C_k^\top \mathbb{E}[y^*_{2,k+1} \mid \mathcal{F}_{k-1}] + F_k^\top \mathbb{E}[y^*_{2,k+1} \omega_k \mid \mathcal{F}_{k-1}] \\
					&\quad + (M_k+S_k\Pi^{2*}_{k}) x^{*}_k + S_k \Sigma^{2}_k + \theta_k = 0,
				\end{split}
			\end{equation}
			where $\mathbf{x}^{*}$ solves the closed-loop system:
			\begin{equation}\label{eq:40.6}
				\begin{cases}
					\begin{split}
						x^*_{k+1} &= \big(A_k+B_k\Pi^{1*}_{k}+C_k\Pi^{2*}_{k}\big)x^*_k + B_k \Sigma^{1*}_{k} \\
						&\quad+ C_k \Sigma^{2*}_{k}  + b_k + \big[\big(D_k+E_k\Pi^{1*}_{k}+F_k\Pi^{2*}_{k}\big)x^*_k \\
						&\quad+ E_k \Sigma^{1*}_{k}  + F_k \Sigma^{2*}_{k} + \sigma_k\big]\omega_k \\
						x^*_0 &= \xi.
					\end{split}
				\end{cases}
			\end{equation}
			
			\noindent\textbf{Analysis of Player 1's Dynamic Behavior}
			
			For any fixed strategy $\left(\Pi^{2*}, \Sigma^{2*}\right) \in \mathscr{M}^{2}(\mathcal{T})$ and initial state $\xi \in \mathbb{R}^{n}$, the FBS$\Delta$E governing Player 1's optimal control admits an adapted solution $(\Sigma^{1*},\mathbf{x}^*,y^*_1)$ satisfying:
			\begin{equation}\label{eq:4.12}
				\left\{
				\begin{aligned}
					x^*_{k+1} &=\big(A_k+B_k\Pi^{1*}_{k}+C_k\Pi^{2*}_{k}\big)x^*_k + B_k \Sigma^{1*}_{k} \\
					&\quad+ C_k \Sigma^{2*}_{k}  + b_k + \big[\big(D_k+E_k\Pi^{1*}_{k}+F_k\Pi^{2*}_{k}\big)x^*_k \\
					&\quad+ E_k \Sigma^{1*}_{k}  + F_k \Sigma^{2*}_{k} + \sigma_{k}\big]\omega_{k}, \\
					y^*_{1,k} &= \left(A_k + C_k \Pi^{2*}_k\right)^\top \mathbb{E}\left[y^*_{1,k+1} \mid \mathcal{F}_{k-1}\right] \\
					&\quad + \left(D_k + F_k \Pi^{2*}_k\right)^\top \mathbb{E}\left[y^*_{1,k+1}\omega_k \mid \mathcal{F}_{k-1}\right] \\
					&\quad + \big(Q_k+L^\top_k\Pi^{1*}_{k}\big) x^*_k +  L^{\top}_k\Sigma^{1*}_{k} + q_k, \\
					x^*_0 &= \xi, \quad y^*_{1,N} = G_N x^*_N + g_N, \quad k \in \mathcal{T},
				\end{aligned}
				\right.
			\end{equation}
			with the stationarity condition
			\begin{equation}\label{eq:4.13}
				\begin{split}
					&B_{k}^{\top}\mathbb{E}[y^*_{1, k+1} \mid \mathcal{F}_{k-1}] + E_{k}^{\top}\mathbb{E}[y^*_{1, k+1}\omega_{k} \mid \mathcal{F}_{k-1}] \\
					&\quad + \big(L_{k}+R_{k}\Pi^{1*}_{k}\big)x^*_{k} + R_{k}\Sigma^{1*}_{k}+\rho_{k} = 0, \quad k\in \mathcal{T}.
				\end{split}
			\end{equation}
			
			For any $\xi\in \mathbb{R}^{n}$, FBS$\Delta$E \eqref{eq:4.12} admits a solution and $(\Pi^{1*},\Sigma^{1*})$ is $\xi$-independent, by subtracting solutions corresponding to $\xi$ and 0, the latter from the former , the adapted solution $(\mathbb{X}^*,\mathbb{Y}^*_1)$  to the following FBS$\Delta$E holds for any $\xi\in \mathbb{R}^{n}$:
			\begin{equation}\label{eq:4.14}
				\left\{
				\begin{aligned}
					\mathbb{X}^*_{k+1} &= \big(A_{k}+B_{k}\Pi^{1*}_{k}+C_{k}\Pi^{2*}_{k}\big)\mathbb{X}^*_{k}\\
					&\quad + \big[\big(D_{k}+E_{k}\Pi^{1*}_{k}+F_{k}\Pi^{2*}_{k}\big)\mathbb{X}^*_{k}\big]\omega_{k}, \\
					\mathbb{Y}^*_{1, k} &= \big(A_{k}+C_{k}\Pi^{2*}_{k}\big)^{\top}\mathbb{E}[\mathbb{Y}^*_{1, k+1} \mid \mathcal{F}_{k-1}] \\
					&\quad + \big(D_{k}+F_{k}\Pi^{2*}_{k}\big)^{\top}\mathbb{E}[\mathbb{Y}^*_{1, k+1}\omega_{k} \mid \mathcal{F}_{k-1}] \\
					&\quad + \big(Q_{k}+L_{k}^{\top}\Pi^{1*}_{k}\big)\mathbb{X}^*_{k}, \\
					\mathbb{X}^*_{0} &= I, \quad \mathbb{Y}^*_{1, N}=G_{N}\mathbb{X}^*_{N}, \quad k\in \mathcal{T},
				\end{aligned}
				\right.
			\end{equation}
		  satisfies
			\begin{equation}\label{eq:4.15}
				\begin{split}
					&B_{k}^{\top}\mathbb{E}[\mathbb{Y}^*_{1, k+1} \mid \mathcal{F}_{k-1}] + E_{k}^{\top}\mathbb{E}[\mathbb{Y}^*_{1, k+1}\omega_{k} \mid \mathcal{F}_{k-1}] \\
					&\quad + \big(L_{k}+R_{k}\Pi^{1*}_{k}\big)\mathbb{X}^*_{k}  = 0, \quad k\in \mathcal{T}.
				\end{split}
			\end{equation}
		
			Then we obatin the desired result \eqref{eq:4.2}-\eqref{eq:4.3}.
			
			\noindent\textbf{Analysis of Player 2's Dynamic Behavior}
			
			Similarly, for any fixed strategy $\left(\Pi^{1*}, \Sigma^{1*}\right) \in \mathscr{M}^{1}(\mathcal{T})$ and initial state $\xi \in \mathbb{R}^{n}$, the adapted solution $(\mathbb{X}^*,\mathbb{Y}^*_2)$  to the following FBS$\Delta$E holds for any $\xi\in \mathbb{R}^{n}$:
			\begin{equation}\label{eq:4.110}
				\left\{
				\begin{aligned}
					\mathbb{X}^*_{k+1} &= \left(A_k + B_k \Pi^{1*}_k+C_{k}\Pi^{2*}_{k}\right)\mathbb{X}^*_{k} \\
					&\quad +\big[\big(D_{k}+E_{k}\Pi^{1*}_{k}+F_{k}\Pi^{2*}_{k}\big)\mathbb{X}^*_{k}\big]\omega_{k}, \\
					\mathbb{Y}^*_{2, k} &= \big(A_{k}+B_{k}\Pi^{1*}_{k}\big)^{\top}\mathbb{E}[\mathbb{Y}^*_{2, k+1} \mid \mathcal{F}_{k-1}] \\
					&\quad + \big(D_{k}+E_{k}\Pi^{1*}_{k}\big)^{\top}\mathbb{E}[\mathbb{Y}^*_{2, k+1}\omega_{k} \mid \mathcal{F}_{k-1}] \\
					&\quad + \big(P_{k}+M_{k}^{\top}\Pi^{2*}_{k}\big)\mathbb{X}^*_{k}, \\
					\mathbb{X}^*_{0} &= I, \quad \mathbb{Y}^*_{2, N}=H_{N}\mathbb{X}^*_{N}, \quad k\in \mathcal{T},
				\end{aligned}
				\right.
			\end{equation}
			satisfies
			\begin{equation}\label{eq:4.120}
				\begin{split}
					&C_{k}^{\top}\mathbb{E}[\mathbb{Y}^{\ast}_{2, k+1} \mid \mathcal{F}_{k-1}] + F_{k}^{\top}\mathbb{E}[\mathbb{Y}^{\ast}_{2, k+1}\omega_{k} \mid \mathcal{F}_{k-1}] \\
					&\quad + \big(M_{k}+S_k\Pi^{2*}_{k}\big)\mathbb{X}^{\ast}_{k} = 0, \quad k\in \mathcal{T}.
				\end{split}
			\end{equation}
		
			Then we obatin the desired result \eqref{eq:4.4}-\eqref{eq:4.5}.~The proof is complete.
		\end{proof}
		
		Next, two Lyapunov-type equations for each Player via dynamic programming principle approach (DPP, for short) is introduced.
		
		For Player 1, given any time step $k\in\mathcal{T}$ and Player 2's optimal strategy $\mathbf{v}^{*}=\Pi^{2*}\mathbf{x}^{*}+\Sigma^{2*}$ (closed-loop form), let $\mathbf{u}^{*}$ be the corresponding optimal control and $y^*_{1}$ its associated dual variable satisfying the following condition:
		\begin{equation}\label{eq:4.18}
			\begin{split}
				&B_{k}^{\top}\mathbb{E}[y^*_{1, k+1} \mid \mathcal{F}_{k-1}] + E_{k}^{\top}\mathbb{E}[y^*_{1, k+1}\omega_{k} \mid \mathcal{F}_{k-1}] \\
				&+ L_{k}x_{k}^{*} + R_{k}u_{k}^{*}+\rho_{k} = 0,\qquad k\in\mathcal{T}.
			\end{split}
		\end{equation}
		
		Define the Value function satisfying the Bellman Equation as follows:
		\begin{equation}\label{eq:4.19}
			\begin{split}
				V_1&(k,x_k^{*}) \\
				=& \frac{1}{2} \mathbb{E}\bigg\{\langle Q_kx_k^{*},x_k^{*}\rangle + 2\langle L_k^\top u_k^{*},x_k^{*}\rangle + \langle R_ku_k^{*},u_k^{*}\rangle \\
				&+ 2\left\langle q_k, x_k^{*}\right\rangle+2\left\langle\rho_k, u_k^{*}\right\rangle + \mathbb{E}\big[ V_1(k+1,x_{k+1}^{*}) \Big| \mathcal{F}_{k-1} \big] \bigg\}.
			\end{split}
		\end{equation}
		
		Assuming a quadratic form for the Value function:
		\begin{equation}\label{eq:4.20}
			V_1(k,x_k^{*}) = \frac{1}{2} \mathbb{E}\bigg[\langle T_k^1 x_k^{*}, x_k^{*} \rangle + 2\langle \phi_k^1, x_k^{*} \rangle + \eta_k^1 \Big| \mathcal{F}_{k-1}\bigg],\quad k\in\mathcal{T},
		\end{equation}
		with terminal conditions:
		\begin{equation}\label{eq:5.12}
			T_N^1 = G_N, \quad \phi_N^1 = g_N, \quad \eta_N^1 = 0.
		\end{equation}
		
		In this Markovian system where the Value function $V_{1}$ is smooth. By \eqref{eq:4.20}, the consistency between the minimization conditions of SMP and DPP implies that:
		\begin{equation}\label{eq:4.220}
			\begin{split}
				y^*_{1,k} &= \nabla_{x_k^{*}} V_{1}(k, x_{k}^{*}) \\
				&= T^{1}_{k}x_{k}^{*}+\phi^{1}_{k},
			\end{split}
		\end{equation}
		which is demonstrably consistent with the anasatzes \eqref{eq:2.25} we postulated in section 3.
		
		Given Player 2's strategy $\left(\Pi^{2*}, \Sigma^{2*}\right) \in \mathscr{M}^{2}(\mathcal{T})$, we substitute the system dynamics into the Bellman equation \eqref{eq:4.19}. By aggregating from $k=0$ to $k=N-1$ and using Appendix A, we obtain the following results:
		\begin{equation}\nonumber
			\begin{split}
				V_1&(k,x_k^{*})\\=&\frac{1}{2}\mathbb{E}\biggl\{\bigg[ \langle Q_kx_k^{*},x_k^{*}\rangle + 2\langle L_k^\top u_k^{*},x_k^{*}\rangle + \langle R_ku_k^{*},u_k^{*}\rangle \\
				&+ 2\left\langle q_k, x_k^{*}\right\rangle +2\left\langle\rho_k, u_k^{*}\right\rangle\bigg]   + \mathbb{E}\Bigg[\langle T_{k+1}^1 x_{k+1}^{*}, x_{k+1}^{*} \rangle   \\
				&+ 2\langle \phi_{k+1}^1, x_{k+1}^{*} \rangle + \eta_{k+1} \Big| \mathcal{F}_{k-1} \Bigg] \biggr\}
				\\=&\frac{1}{2}\mathbb{E}\biggl\{\bigg[ \langle Q_kx_k^{*},x_k^{*}\rangle + 2\langle L_k^\top u_k^{*},x_k^{*}\rangle + \langle R_ku_k^{*},u_k^{*}\rangle \\
				&+ 2\left\langle q_k, x_k^{*}\right\rangle +2\left\langle\rho_k, u_k^{*}\right\rangle\bigg]   + \mathbb{E}\Bigg[\langle T_{k+1}^1 \big(A_{k} x_{k}^{*}+B_{k} u_{k}^{*}\\&+C_{k} v^{*}_{k}+b_{k}+\left(D_{k} x_{k}^{*}+E_{k} u_{k}^{*}+F_{k} v^{*}_{k}+\sigma_{k}\right) \omega_{k}\big), A_{k} x_{k}^{*}\\&+B_{k} u_{k}^{*}+C_{k} v^{*}_{k}+b_{k}+\left(D_{k} x_{k}^{*}+E_{k} u_{k}^{*}+F_{k} v^{*}_{k}+\sigma_{k}\right) \omega_{k} \rangle   \\&+ 2\langle \phi_{k+1}^1, A_{k} x_{k}^{*}+B_{k} u_{k}^{*}+C_{k} v^{*}_{k}+b_{k}+\big(D_{k} x_{k}^{*}+E_{k} u_{k}^{*}\\&+F_{k} v^{*}_{k}+\sigma_{k}\big) \omega_{k} \rangle + \eta_{k+1} \Big| \mathcal{F}_{k-1} \Bigg] \biggr\}
			\\=&\frac{1}{2}\mathbb{E}\biggl\{\bigg[ \langle Q_kx_k^{*},x_k^{*}\rangle + 2\langle L_k^\top u_k^{*},x_k^{*}\rangle + \langle R_ku_k^{*},u_k^{*}\rangle \\&+ 2\left\langle q_k, x_k^{*}\right\rangle +2\left\langle\rho_k, u_k^{*}\right\rangle\bigg]   + \mathbb{E}\Bigg[\langle T_{k+1}^1 \big(A_{k} x_{k}^{*}+B_{k} u_{k}^{*}\\&+C_{k} \big(\Pi^{2*}_{k}x_{k}^{*}+\Sigma^{2*}_{k}\big)+b_{k}+\big(D_{k} x_{k}^{*}+E_{k} u_{k}^{*}+F_{k}	\\&\times \big(\Pi^{2*}_{k}x_{k}^{*}+\Sigma^{2*}_{k}\big)+\sigma_{k}\big) \omega_{k}\big), A_{k} x_{k}^{*}+B_{k} u_{k}^{*}+C_{k}\\&\times \big(\Pi^{2*}_{k}x_{k}^{*}+\Sigma^{2*}_{k}\big)+b_{k}+\big(D_{k} x_{k}^{*}+E_{k} u_{k}^{*}+F_{k} \big(\Pi^{2*}_{k}\\&\times x_{k}^{*}+\Sigma^{2*}_{k}\big)+\sigma_{k}\big) \omega_{k} \rangle   + 2\langle \phi_{k+1}^1, A_{k} x_{k}^{*}+B_{k} u_{k}^{*}+C_{k}\\&\times \big(\Pi^{2*}_{k}x_{k}^{*}+\Sigma^{2*}_{k}\big)+b_{k}+\big(D_{k} x_{k}^{*}+E_{k} u_{k}^{*}+F_{k} \big(\Pi^{2*}_{k}x_{k}^{*}\\& +\Sigma^{2*}_{k}\big)+\sigma_{k}\big)\omega_{k} \rangle + \eta_{k+1} \Big| \mathcal{F}_{k-1} \Bigg]\biggr\}\\
	=&\frac{1}{2}\mathbb{E}\biggl\{\bigg[ 2\langle L_k^\top u_k^{*},x_k^{*}\rangle + \langle R_ku_k^{*},u_k^{*}\rangle + 2\left\langle q_k, x_k^{*}\right\rangle \\
				&+2\left\langle\rho_k, u_k^{*}\right\rangle\bigg]   + \mathbb{E}\Bigg[\big\langle \Delta\left(T^{1}_{k+1},\Pi^{2*}_{k}\right)x_k^{*},x_k^{*}\big\rangle + \big\langle T_{k+1}^1\\&\times\big[ A_{k}+C_k\Pi^{2*}_{k}+\big(D_{k}+F_k\Pi^{2*}_{k}\big)\omega_{k}\big]x_{k}^{*},B_{k} u_{k}^{*}+C_k\\&\times\Sigma^{2*}_{k}+b_k+\big(E_ku_k^{*}+F_k\Sigma^{2*}_{k}+\sigma_{k}\big)\omega_{k}\big\rangle+\big\langle T_{k+1}^1
\\
&\times\big[B_{k} u_{k}^{*}+C_k\Sigma^{2*}_{k} +b_k+\big(E_ku_{k}^{*}+F_k\Sigma^{2*}_{k}+\sigma_{k}\big)\omega_{k}\big],\\&\big(A_{k} x_{k}^{*}+B_{k} u_{k}^{*}+C_{k} \big(\Pi^{2*}_{k}x_{k}^{*}+\Sigma^{2*}_{k}\big)+b_{k}+\big(D_{k} x_{k}^{*}\end{split}
\end{equation}

\begin{equation}\label{eq:5.13}
\begin{split}&+E_{k} u_{k}^{*}+F_{k} \big(\Pi^{2*}_{k}x_{k}^{*}+\Sigma^{2*}_{k}\big)+\sigma_{k}\big) \omega_{k}\big)\big\rangle + 2\big\langle \phi_{k+1}^1, \\&  A_{k} x_{k}^{*}+B_{k} u_{k}^{*}+C_{k} \big(\Pi^{2*}_{k}x_{k}^{*}+\Sigma^{2*}_{k}\big)+b_{k}+\big(D_{k} x_{k}^{*}+\\&E_{k} u_{k}^{*}+F_{k} \big(\Pi^{2*}_{k}x_{k}^{*}+\Sigma^{2*}_{k}\big)+\sigma_{k}\big) \omega_{k} \big\rangle + \eta_{k+1} \Big| \mathcal{F}_{k-1} \Bigg]\biggr\}.
			\end{split}
		\end{equation}
		Substituting the closed-loop control law $\mathbf{u}^{*}=\Pi^{1*}\mathbf{x}^{*}+\Sigma^{1*}$ into \eqref{eq:5.13}, then we have
		\begin{equation}\nonumber
			\begin{split}
				&V_1(k,x_k^{*})\\
				=&\frac{1}{2}\mathbb{E}\biggl\{\bigg[ 2\langle L_k^\top (\Pi^{1*}_{k}x_k^{*}+\Sigma^{1*}_{k}),x_k^{*}\rangle + \langle R_k(\Pi^{1*}_{k}x_k^{*}+\Sigma^{1*}_{k}),\\
				&(\Pi^{1*}_{k}x_k^{*}+\Sigma^{1*}_{k})\rangle + 2\left\langle q_k, x_k^{*}\right\rangle +2\left\langle\rho_k, \Pi^{1*}_{k}x_k^{*}+\Sigma^{1*}_{k}\right\rangle\bigg]   \\
				&+ \mathbb{E}\Bigg[ \big\langle \Delta\left(T^{1}_{k+1},\Pi^{2*}_{k}\right)x_k^{*},x_k^{*}\big\rangle
				+ \big\langle T_{k+1}^1  \big[ A_{k} + C_k\Pi^{2*}_{k} \\
				&+ \big(D_{k} + F_k\Pi^{2*}_{k}\big)\omega_{k}\big]x_{k}^{*},
				B_{k} (\Pi^{1*}_{k}x_k^{*} + \Sigma^{1*}_{k}) + C_k \Sigma^{2*}_{k} 	\\
				&+ b_k + \big(E_k(\Pi^{1*}_{k}x_k^{*} + \Sigma^{1*}_{k}) + F_k\Sigma^{2*}_{k} + \sigma_{k}\big)\omega_{k} \big\rangle
				+ \big\langle T_{k+1}^1\\
				&\times \big[ B_{k} (\Pi^{1*}_{k} x_k^{*} + \Sigma^{1*}_{k})	 + C_k\Sigma^{2*}_{k} + b_k
				+ \big(E_k (\Pi^{1*}_{k}x_k^{*} + \Sigma^{1*}_{k})\\
				&+ F_k\Sigma^{2*}_{k} + \sigma_{k}\big)\omega_{k}\big],
				\big( A_{k} x_{k}^{*} + B_{k} (\Pi^{1*}_{k}x_k^{*} + \Sigma^{1*}_{k}) + C_{k}\\
				&\times \big(\Pi^{2*}_{k}x_{k}^{*} + \Sigma^{2*}_{k}\big) + b_{k}
				+ \big( D_{k} x_{k}^{*} + E_{k} (\Pi^{1*}_{k}x_k^{*} + \Sigma^{1*}_{k})\\
				&+ F_{k} \big(\Pi^{2*}_{k}x_{k}^{*} + \Sigma^{2*}_{k}\big) + \sigma_{k}\big) \omega_{k} \big) \big\rangle
				+ 2\big\langle \phi_{k+1}^1, A_{k} x_{k}^{*} + B_{k}\\
				&\times (\Pi^{1*}_{k}x_k^{*} + \Sigma^{1*}_{k}) + C_{k} \big(\Pi^{2*}_{k}x_{k}^{*}  + \Sigma^{2*}_{k}\big) + b_{k}
				+ \big( D_{k} x_{k}^{*}  \\
				&+ E_{k} (\Pi^{1*}_{k}x_k^{*} + \Sigma^{1*}_{k}) + F_{k} \big(\Pi^{2*}_{k}x_{k}^{*} +\Sigma^{2*}_{k}\big) + \sigma_{k}\big) \omega_{k} \big\rangle
				\\&+ \eta_{k+1} \Big| \mathcal{F}_{k-1} \Bigg]\biggr\}
					\end{split}
			\end{equation}
		\begin{equation}\label{eq:4.23}
	\begin{split}
				=&\frac{1}{2}\mathbb{E}\biggl\{\mathbb{E}\bigg[ \bigg\langle\bigg( \Delta\left(T^{1}_{k+1},\Pi^{2*}_{k}\right)
				+\Pi^{1*\top}_{k}\Upsilon(T^{1}_{k+1})\Pi^{1*}_{k}	\\
				&+\Lambda\left(T^{1}_{k+1},\Pi^{2*}_{k}\right)^{\top}\Pi^{1*}_{k}+\Pi^{1*\top}_{k}\Lambda\left(T^{1}_{k+1},\Pi^{2*}_{k}\right)\bigg)x_k^{*},x_k^{*}\big\rangle\\
				&+2\bigg\langle \Theta\left(T^{1}_{k+1},\Pi^{2*}_{k},\phi^{1}_{k+1}\right)
				+\Lambda^{\top}(T^{1}_{k+1},\Pi^{2*}_{k})	\Sigma^{1*}_{k}	+	\Pi^{1*\top}_{k}\\
				&\times\big[\Upsilon\left(T_{k+1}^{1}\right)\Sigma^{1*}_{k}+\Phi\left(T_{k+1}^{1}, \Pi_{k}^{2*}, \phi_{k+1}^{1}\right)\big],x_k^{*}\big\rangle
		\\
				&+\big\langle R_k\Sigma^{1*}_{k} ,\Sigma^{1*}_{k}\big\rangle+2\big\langle \rho_k ,\Sigma^{1*}_{k}\big\rangle + \bigg\langle T^{\top}_{k+1}\big[B_k\Sigma^{2*}_{k}+C_k\Sigma^{2*}_{k}\\
				&+b_k+(E_k+F_k\Sigma^{2*}_{k}+\sigma_{k})\omega_{k}\big] ,B_k\Sigma^{1*}_{k}+C_k\Sigma^{2*}_{k}+b_k\\
				&+E_k\Sigma^{1*}_{k}\omega_{k}+F_k\Sigma^{2*}_{k}\omega_{k}+b_k\omega_{k}\bigg\rangle  \Big| \mathcal{F}_{k-1} \Bigg]\biggr\}.
			\end{split}
		\end{equation}
		Comparing the quadratic terms in $x_k$ between \eqref{eq:4.20} and \eqref{eq:4.23}, for any $k\in\mathcal{T}$, the cross-coupled Lyapunov-type equations for $T^1$ can be derived as:
		\begin{equation}\label{eq:50.14}
			\left\{
			\begin{aligned}
				T^{1}_{k} &= \Delta\left(T^{1}_{k+1},\Pi^{2*}_{k}\right)
				+ \Pi^{1\top}_{k}\Upsilon(T^{1}_{k+1})\Pi^{1*}_{k} \\
				&\quad + \Lambda\left(T^{1}_{k+1},\Pi^{2*}_{k}\right)^{\top}\Pi^{1*}_{k}
				+ \Pi^{1*\top}_{k}\Lambda\left(T^{1}_{k+1},\Pi^{2*}_{k}\right), \\
				T^{1}_{N} &= G_{N}.
			\end{aligned}
			\right.
		\end{equation}
	Similarly, cross-coupled Lyapunov-type equations for $T^2$  can be derived as:
		\begin{equation}\label{eq:50.15}
			\left\{
			\begin{aligned}
				T^{2}_{k} &= \Delta\left(T^{2}_{k+1},\Pi^{1*}_{k}\right)
				+ \Pi^{2*\top}_{k}\Upsilon(T^{2}_{k+1})\Pi^{2}_{k} \\
				&\quad + \Lambda\left(T^{2}_{k+1},\Pi^{2*}_{k}\right)^{\top}\Pi^{2*}_{k}
				+ \Pi^{2*\top}_{k}\Lambda\left(T^{2}_{k+1},\Pi^{1*}_{k}\right), \\
				T^{2}_{N} &= H_{N}.
			\end{aligned}
			\right.
		\end{equation}
		
		Substituting the relationship \eqref{eq:4.220} and closed-loop forms \eqref{eq:2.9} of $\mathbf{u}^{*}$ and $\mathbf{v}^*$ into \eqref{eq:4.18}, the we can deduce the following equation by using notations in Appendix A:
		\[
		\Upsilon\left(T_{k+1}^{1}\right)\Sigma^{1*}_{k} + \Phi\left(T_{k+1}^{1}, \Pi_{k}^{2*}, \phi_{k+1}^{1}\right) = 0, \quad \forall k \in \mathcal{T}.
		\]
		
		Furthermore, comparing the linear terms in $x_k$ gives BS$\Delta$Es for $\phi^1$:
		\begin{equation}\label{eq:4.190}
			\left\{
			\begin{aligned}
				\phi^{1}_{k} &= \Theta\left(T^{1}_{k+1},\Pi^{2*}_{k},\phi^{1}_{k+1}\right)
				+ \Lambda^{\top}\left(T^{1}_{k+1},\Pi^{2*}_{k}\right)\Sigma^{1*}_{k}, \\
				\phi^{1}_{N} &= g_{N}.
			\end{aligned}
			\right.
		\end{equation}
	Similarly, BS$\Delta$Es for $\phi^2$ can be derived as
		\begin{equation}\label{eq:4.200}
			\left\{
			\begin{aligned}
				\phi^{2}_{k} &= \Theta\left(T^{2}_{k+1},\Pi^{1*}_{k},\phi^{2}_{k+1}\right)
				+ \Lambda^{\top}\left(T^{2}_{k+1},\Pi^{1*}_{k}\right)\Sigma^{2*}_{k}, \\
				\phi^{2}_{N} &= h_{N}.
			\end{aligned}
			\right.
		\end{equation}
		
		Next, the regular solvability of the closed-loop Nash equilibrium for DSDG-CL boils down to another main result.
		
		\begin{thm}\label{thm:5.2}
			Let Assumptions \ref{ass:2.1}-\ref{ass:2.2} hold.~Strategy quadruple \(\left(\Pi^{1*}, \Sigma^{1*}; \Pi^{2*},  \Sigma^{2*}\right)\) \(\in \mathscr{M}(\mathcal{T})\) is a closed-loop Nash equilibrium of Problem (DSDG-CL) on $\mathcal{T}$ if and only if the following hold:
			
			(i) For Player1, the solution $T^{1} \in  L_{\mathbb{F}}^{\infty}(\mathcal{\overline{T}}; \mathbb{S}_{+}^{n})$ to the Lyapunov-type equation
			\begin{equation}\label{eq:5.15}
				\left\{\begin{array}{lll}
					\begin{split}
						T^{1}_{k}=&	\Delta\left(T^{1}_{k+1},\Pi^{2\ast}_{k}\right)
						+
						\Pi^{1\ast\top}_{k}\Upsilon(T^{1}_{k+1}) \Pi^{1\ast}_{k}+\\&\Lambda^{\top}\left(T^{1}_{k+1},\Pi^{2\ast}_{k}\right)\Pi^{1\ast}_{k}+\Pi^{1\ast\top}_{k}\Lambda\left(T^{1}_{k+1},\Pi^{2\ast}_{k}\right),
						\\
						T^{1}_{N}=&G_{N}.
					\end{split}
				\end{array}
				\right.
			\end{equation}
			satisfies the following two conditions:
			\begin{equation}\label{eq:5.16}
				\begin{split}
					&\Upsilon(T^{1}_{k+1}) \succeq 0.
				\end{split}
			\end{equation}
			\begin{equation}\label{eq:5.17}
				\begin{split}
					&\Lambda\left(T^{1}_{k+1},\Pi^{2\ast}_{k}\right)+\Upsilon(T^{1}_{k+1}) \Pi^{1\ast}_{k}=0,\quad k\in\mathcal{T}.
				\end{split}
			\end{equation}
			And the adapted solution $\phi^{1}$ to the BS$\Delta$Es
		{\small\begin{equation}\label{eq:5.18}
					\left\{
					\begin{aligned}
						\begin{split}
							\phi_{k}^{1} =& \Lambda^{\top}\left(T_{k+1}^{1}, \Pi_{k}^{2*}\right)\Sigma^{1*}_{k}+\Theta\left(T_{k+1}^{1}, \Pi_{k}^{2*}, \phi_{k+1}^{1}\right),
							\\
							\phi_{N}^{1} =& g_{N}.
						\end{split}
					\end{aligned}
					\right.
			\end{equation}}
			satisfies
			\begin{equation}\label{eq:5.19}
				\begin{split}
					&\Upsilon\left(T_{k+1}^{1}\right)\Sigma^{1*}_{k}+\Phi\left(T_{k+1}^{1}, \Pi_{k}^{2*}, \phi_{k+1}^{1}\right)
					= 0,\quad k\in\mathcal{T}.
				\end{split}
			\end{equation}
			
			(ii) For Player2, the solution $ T^{2} \in  L_{\mathbb{F}}^{\infty}
			(\overline{\mathcal{T}}; \mathbb{S}_{+}^{n}) $ to the Lyapunov-type equation
			\begin{equation}\label{eq:5.20}
				\left\{\begin{array}{lll}
					\begin{split}
						T^{2}_{k}
						=&\Delta\left(T^{2}_{k+1},\Pi^{1*}_{k}\right)+\Pi^{2*\top}_{k}\Upsilon(T^{2}_{k+1})\Pi^{2*}_{k}+\\&\Lambda^{\top}\left(T^{2}_{k+1},\Pi^{1*}_{k}\right)\Pi^{2*}_{k}+\Pi^{2*\top}_{k}\Lambda\left(T^{2}_{k+1},\Pi^{1*}_{k}\right),
						\\
						T^{2}_{N}=&H_{N}.
					\end{split}
				\end{array}
				\right.
			\end{equation}
			satisfies the following two conditions:
			\begin{equation}\label{eq:5.21}
				\begin{split}
					&\Upsilon(T^{2}_{k+1})\succeq 0.
				\end{split}
			\end{equation}
			\begin{equation}\label{eq:5.22}
				\begin{split}
					&\Lambda\left(T^{2}_{k+1},\Pi^{1\ast}_{k}\right)+\Upsilon(T^{2}_{k+1})\Pi^{2*}_{k}=0,\quad k\in\mathcal{T}.
				\end{split}
			\end{equation}
			And the adapted solution $\phi^{1}$ to the BS$\Delta$Es
			\begin{equation}\label{eq:5.23}
				\left\{\begin{array}{lll}
					\begin{split}
						\phi_{k}^{2}=&  \Lambda^{\top}\left(T_{k+1}^{2}, \Pi_{k}^{1*}\right)\Sigma^{2*}_{k}+\Theta\left(T_{k+1}^{2}, \Pi_{k}^{1*}, \phi_{k+1}^{2}\right),
						\\
						\phi_{N}^{2}=&h_{N}.
					\end{split}
				\end{array}
				\right.
			\end{equation}
			satisfies
			\begin{equation}\label{eq:5.24}
				\begin{split}
					& \Upsilon\left(T_{k+1}^{2}\right)\Sigma^{2*}_{k}+\Phi\left(T_{k+1}^{2}, \Pi_{k}^{1*}, \phi_{k+1}^{2}\right)
					= 0,\quad k\in\mathcal{T}.
				\end{split}
			\end{equation}
		\end{thm}
		\begin{proof}
			
			We first prove the necessity. For Player1, suppose that \(\left(\Pi^{1*}, \Sigma^{1*}; \Pi^{2*},  \Sigma^{2*}\right)\) is a closed-loop Nash equilibrium of DSDG-CL.~Based on the Proposition \ref{prop4.1}, we can obtain the FBS$\Delta$E \eqref{eq:4.14} satisfying \eqref{eq:4.15}.
			 According to the discussion in Section 3, equation \eqref{eq:4.14} can be decoupled through the Lyapunov-type equation \eqref{eq:5.15}. i.e.
			
			\[	\mathbb{Y}^*_{1, k}=T_{k}^{1} 	\mathbb{X}^*_{k}.\]
		
			Substituting above into \eqref{eq:4.15} and using notations in Appendix A, we obtain:
			\begin{equation}\label{eq:5.28}
				\begin{split}
					0= & B_{k}^{\top}\mathbb{E}[\mathbb{Y}^*_{1, k+1} \mid \mathcal{F}_{k-1}] + E_{k}^{\top}\mathbb{E}[\mathbb{Y}^*_{1, k+1}\omega_{k} \mid \mathcal{F}_{k-1}] \\
					& + \big(L_{k}+R_{k}\Pi^{1*}_{k}\big)\mathbb{X}^*_{k}   \\
					= & B_{k}^{\top} \mathbb{E}\left[T_{k+1}^{1} 	\mathbb{X}^*_{k+1} \mid \mathcal{F}_{k-1}\right]+E_{k}^{\top} \mathbb{E}\big[\left(T_{k+1}^{1} 	\mathbb{X}^*_{k+1}\right) \\&\times\omega_{k} \mid \mathcal{F}_{k-1}\big] + \big(L_{k}+R_{k}\Pi^{1*}_{k}\big)\mathbb{X}^*_{k}
					\\
					=& B_{k}^{\top} \mathbb{E}\bigg[T_{k+1}^{1}\bigg(\big(A_{k}+C_{k}\Pi^{2*}_{k}\big)	\mathbb{X}^*_{k}+B_{k}\Pi^{1*}_{k}	\mathbb{X}^*_{k}+\\&\bigg[\big(D_{k}+F_{k}\Pi^{2*}_{k}\big) 	\mathbb{X}^*_{k}+E_{k}\Pi^{1*}_{k}	\mathbb{X}^*_{k}\bigg]\omega_{k}\bigg) \mid \mathcal{F}_{k-1}\bigg]\\&+E_{k}^{\top} \mathbb{E}\bigg[\bigg(T_{k+1}^{1}\bigg(\big(A_{k}+C_{k}\Pi^{2*}_{k}\big)	\mathbb{X}^*_{k}+B_{k}\Pi^{1*}_{k}\mathbb{X}^*_{k}	\\&+\bigg[\big(D_{k}+F_{k}\Pi^{2*}_{k}\big)\mathbb{X}^*_{k}+E_{k}\Pi^{1*}_{k}\mathbb{X}^*_{k}\bigg]\omega_{k}\bigg)\bigg) \omega_{k} \mid \mathcal{F}_{k-1}\bigg]\\&+\big(L_{k}+R_{k}\Pi^{1*}_{k}\big)\mathbb{X}^*_{k}
					\\
					= & \big[\Lambda\left(T^{1}_{k+1},\Pi^{2\ast}_{k}\right)+\Upsilon(T^{1}_{k+1}) \Pi^{1\ast}_{k}\big]\mathbb{X}_{k}.
				\end{split}
			\end{equation}
			
			The almost sure invertibility of $\mathbb{X}^*_k$ implies that for any realization $\omega \in \Omega$ where $\mathbb{X}^*_k(\omega) \neq 0$, the coefficient matrix must satisfy $\eqref{eq:5.17}$.
			
			Now we suppose $\phi^{1}$ is the adapted solution to BS$\Delta$Es \eqref{eq:5.18} and let $(\mathbf{x}^*,y^*_{1,k})$ be the solutuion to \eqref{eq:2.22}.
			
			In the following, we can prove
		\begin{equation}\label{eq:relationship}
			\phi_{k}^{1}=y^*_{1, k}-T_{k}^{1} x^*_{k}.
		\end{equation}				
			
			In fact, for any $\Pi^{1}\in\mathscr{M}^{1}(\mathcal{T})$, $\phi_{N}^{1}=g_{N}$ , we have
			\small\begin{equation}\label{eq:5.30}
				\begin{split}
					&  y^*_{1, k} - T_{k}^{1} x^*_{k}\\=& \big(A_{k}+C_{k}\Pi^{2*}_{k}\big)^{\top}\mathbb{E}[y^*_{1, k+1} \mid \mathcal{F}_{k-1}] + \big(D_{k}+F_{k}\Pi^{2*}_{k}\big)^{\top}\\&\mathbb{E}[y^*_{1, k+1}\omega_{k} \mid \mathcal{F}_{k-1}] + Q_{k}x^*_{k} + L_{k}^{\top}u^*_{k}+q_{k}- T_{k}^{1} x^*_{k}\\=& \big(A_{k}+C_{k}\Pi^{2*}_{k}\big)^{\top}\mathbb{E}[T^{1}_{k+1}x^*_{k+1}+\phi_{k+1} \mid \mathcal{F}_{k-1}] \\&+ \big(D_{k}+F_{k}\Pi^{2*}_{k}\big)^{\top}\mathbb{E}[\big(T^{1}_{k+1}x^*_{k+1}+\phi_{k+1}\big)\omega_{k} \mid \mathcal{F}_{k-1}] \\&+ Q_{k}x^*_{k} + L_{k}^{\top}u^*_{k}+q_{k}- T_{k}^{1} x^*_{k}
				\\
					=& \big(A_{k}+C_{k}\Pi^{2*}_{k}\big)^{\top}\mathbb{E}[T^{1}_{k+1}\big(\big(A_{k}+C_{k}\Pi^{2*}_{k}\big)x^*_{k}\\&+B_{k}u^*_{k}+C_{k}\Sigma^{2*}_{k}+b_{k}+\big[\big(D_{k}+F_{k}\Pi^{2*}_{k}\big)x^*_{k}\\&+E_{k}u^*_{k}+F_{k}\Sigma^{2*}_{k}+\sigma_{k}\big]\omega_{k}\big)+\phi_{k+1} \mid \mathcal{F}_{k-1}] 	\\&+ \big(D_{k}+F_{k}\Pi^{2*}_{k}\big)^{\top}\mathbb{E}[\big(T^{1}_{k+1}\big(\big(A_{k}+C_{k}\Pi^{2*}_{k}\big)x^*_{k}\\&+B_{k}u^*_{k}+C_{k}\Sigma^{2*}_{k}+b_{k}+\big[\big(D_{k}+F_{k}\Pi^{2*}_{k}\big)x^*_{k}\\&+E_{k}u^*_{k}+F_{k}\Sigma^{2*}_{k}+\sigma_{k}\big]\omega_{k}\big)+\phi_{k+1}\big)\omega_{k} \mid \mathcal{F}_{k-1}]	\\& + Q_{k}x^*_{k} + L_{k}^{\top}u^*_{k}+q_{k}- T_{k}^{1} x^*_{k}
			\\=&  \Lambda^{\top}\left(T_{k+1}^{1}, \Pi_{k}^{2*}\right) u^*_{k}+\Delta(T^{1}_{k+1},\Pi^{2*}_{k})x^*_k\\&+\Theta\left(T_{k+1}^{1}, \Pi_{k}^{2*}, \phi_{k+1}^{1}\right)-T^{1}_{k}x^{*}_{k}
						\\=&-\bigg[\Pi^{1*\top}_{k}\Upsilon(T^{1}_{k+1})\Pi^{1*}_{k}+\Pi^{1*\top}_{k}\Lambda(T^{1}_{k+1},\Pi^{2*}_{k})\bigg]x^*_k\\&+\Lambda^{\top}\left(T_{k+1}^{1}, \Pi_{k}^{2*}\right)\Sigma^{1*}_{k}+\Theta\left(T_{k+1}^{1}, \Pi_{k}^{2*}, \phi_{k+1}^{1}\right)
					.
				\end{split}
			\end{equation}
			where the closed-loop control law of $\mathbf{u}^*$ \eqref{eq:2.9} and \eqref{eq:50.14} are employed to prove the final equality.
			
			It is naturally the last equation of \eqref{eq:5.30} is consistent with \eqref{eq:5.18} when \eqref{eq:5.17} satisfied. we've proved the \eqref{eq:relationship}.
			
			From the stationarity condition in  FBS$\Delta$Es~\eqref{eq:2.22},~we have
			\small{\begin{equation}\nonumber
				\begin{split}
					0= & B_{k}^{\top} \mathbb{E}\big[y^*_{1, k+1} \mid \mathcal{F}_{k-1}\big]+E_{k}^{\top} \mathbb{E}\big[y^*_{1, k+1} \omega_{k} \mid \mathcal{F}_{k-1}\big]\\&+L_{k} x^*_{k}+R_{k} u^*_{k}+\rho_{k}
					\\
					= & B_{k}^{\top} \mathbb{E}\big[T_{k+1}^{1} x^*_{k+1}+\phi_{k+1}^{1} \mid \mathcal{F}_{k-1}\big]+E_{k}^{\top} \mathbb{E}\big[\big(T_{k+1}^{1} x^*_{k+1}\\&+\phi_{k+1}^{1}\big) w_{k} \mid \mathcal{F}_{k-1}\big]+L_{k} x^*_{k}+R_{k} u^*_{k}+\rho_{k}
				\\
					= & B_{k}^{\top} \mathbb{E}\big[T_{k+1}^{1}\big[\big(A_{k}+C_{k} \Pi_{k}^{2*}\big) x^*_{k}+B_{k} u^*_{k}+C_{k} \Sigma_{k}^{2*}+b_{k}\\&+\big[\big(D_{k}+F_{k} \Pi_{k}^{2*}\big) x^*_{k}+E_{k} u^*_{k}+F_{k} \Sigma_{k}^{2*}+\sigma_{k}\big] \omega_{k}\big]+\phi_{k+1}^{1} \\
					&\mid \mathcal{F}_{k-1}\big]+E_{k}^{\top} \mathbb{E}\big[\big(T _ { k + 1 } ^ { 1 } \big[\big(A_{k}+C_{k} \Pi_{k}^{2*}\big) x^*_{k} +B_{k} u^*_{k}\\&+C_{k} \Sigma_{k}^{2*}+b_{k}+\big[\big(D_{k}+F_{k} \Pi_{k}^{2*}\big) x^*_{k}+E_{k} u^*_{k}+F_{k} \Sigma_{k}^{2*}\\
					&+\sigma_{k}\big] \omega_{k}\big]+\phi_{k+1}^{1}\big) \omega_{k} \mid \mathcal{F}_{k-1}\big]+L_{k} x^*_{k}+R_{k} u^*_{k}+\rho_{k}\end{split}
				\end{equation}}
			
		\begin{equation}\label{eq:5.32}
	\begin{split}
					=&\Upsilon\left(T_{k+1}^{1}\right) u^*_{k}+\Lambda(T^{1}_{k+1},\Pi^{2*}_{k})x^*_{k}+\Phi\left(T_{k+1}^{1}, \Pi_{k}^{2*}, \phi_{k+1}^{1}\right)
					\\=&\bigg[\Lambda(T^{1}_{k+1},\Pi^{2*}_{k})+\Upsilon\left(T_{k+1}^{1}\right) \Pi^{1*}_{k}\bigg]x^{*}_k+\Upsilon\left(T_{k+1}^{1}\right)\Sigma^{1*}_{k}\\&+\Phi\left(T_{k+1}^{1}, \Pi_{k}^{2*}, \phi_{k+1}^{1}\right),
				\end{split}
			\end{equation}
			where the closed-loop control law of $\mathbf{u}^*$ \eqref{eq:2.9} is employed to prove the final equality.

			Then, we obtain condition \eqref{eq:5.19}, which is consistent with \eqref{eq:5.17}.
			
			Similarly for Player2, we can respectively obtain the \eqref{eq:5.22}, \eqref{eq:5.23} and \eqref{eq:5.24} as well. Besides, the proof of \eqref{eq:5.16} and  \eqref{eq:5.21} will be included in the proof of sufficiency.
			
			To prove sufficiency, let
			$\mathbf{x}=X\big(\Pi^{1*},\Sigma^{1};\Pi^{2*},\Sigma^{2\ast},\xi\big)$ be the state process corresponding to $(0,\xi)$ and $\big(\Pi^{1*},\Sigma^{1};\Pi^{2*},\Sigma^{2\ast},\xi\big)$. Furthermore, the state dynamics of $\mathbf{x}$ are governed by equation \eqref{eq:40.1},
			while the cost functional for Player 1 is given by equation \eqref{eq:40.2}.
			
	     Analyze the cumulative difference in inner products $\langle y_{1,k+1}, x_{k+1}\rangle - \langle y_{1,k}, x_k\rangle$ across steps $k = 0$ to $N-1$ by combining the ansatz from \eqref{eq:2.25} with \eqref{eq:3.13} and using the notations in Appendix A, then we have:
	     \begin{equation}\nonumber
			\begin{split}
				&\frac{1}{2}\mathbb{E}\bigl\{\langle G_{N} x_{N}, x_{N}\rangle + 2\langle g_{N}, x_{N}\rangle\bigr\}-\frac{1}{2}\mathbb{E}\bigl\{\langle T^{1}_{0} x_{0}, x_{0}\rangle\\&+2\langle \phi^{1}_{0} , x_{0}\rangle \bigr\}
				\\=&\frac{1}{2}\mathbb{E}\bigl\{\langle G_{N} x_{N}, x_{N}\rangle - \langle T^{1}_{0}x_{0}, x_{0}\rangle+2\big[\langle g_{N}, x_{N}\rangle-\langle \phi^{1}_{0} , x_{0}\rangle\big]\bigr\}
				\\=&\frac{1}{2} \mathbb{E}\sum_{k=0}^{N-1}\biggl\{ \langle T^{1}_{k+1}x_{k+1},x_{k+1}\rangle-\langle T^{1}_{k}x_{k},x_{k}\rangle+2\big[\langle \phi^{1}_{k+1},x_{k+1}\rangle\\&-\langle \phi^{1}_{k},x_{k}\rangle\big] \biggr\}
					\\=&\frac{1}{2} \mathbb{E}\sum_{k=0}^{N-1}\biggl\{ \langle T^{1}_{k+1}\big(\big(A_k+B_k\Pi^{1*}_{k}+C_k\Pi^{2*}_{k}\big)x_k + B_k \Sigma^{1}_{k} \\&+ C_k \Sigma^{2*}_{k} + b_k + \big[\big(D_k+E_k\Pi^{1*}_{k}+F_k\Pi^{2*}_{k}\big)x_k + E_k \Sigma^{1}_{k} \\&+ F_k \Sigma^{2*}_{k} + \sigma_k\big]\omega_k)\big),\big(A_k+B_k\Pi^{1*}_{k}+C_k\Pi^{2*}_{k}\big) x_k + B_k \Sigma^{1}_{k} \\&+ C_k \Sigma^{2*}_{k} + b_k + \big[\big(D_k+E_k\Pi^{1*}_{k}+F_k\Pi^{2*}_{k}\big)x_k + E_k \Sigma^{1}_{k} \\&+ F_k \Sigma^{2*}_{k} + \sigma_k\big]\omega_k)\rangle-\langle T^{1}_{k}x_{k},x_{k}\rangle+2\langle \phi^{1}_{k+1},\big(A_k+B_k\Pi^{1*}_{k}\\&+C_k\Pi^{2*}_{k}\big) x_k + B_k \Sigma^{1}_{k} + C_k \Sigma^{2*}_{k} + b_k + \big[\big(D_k+E_k\Pi^{1*}_{k}\\&+F_k\Pi^{2*}_{k}\big)x_k + E_k \Sigma^{1}_{k} + F_k \Sigma^{2*}_{k} + \sigma_k\big]\omega_k\rangle-2\langle \phi^{1}_{k},x_{k}\rangle \biggr\}
					\\=&\frac{1}{2} \mathbb{E}\sum_{k=0}^{N-1}\biggl\{\big\langle\big[\Delta(T^{1}_{k+1},\Pi^{2*}_{k})-Q_{k}+\big(\Lambda(T^{1}_{k+1},\Pi^{2*}_{k})\\&\quad-L_{k}\big)^{\top}\Pi^{1*}_{k}+\Pi^{1*\top}_{k}\big(\Lambda(T^{1}_{k+1},\Pi^{2*}_{k})-L_{k}\big)\\&\quad+\Pi^{1*\top}_{k}\times(\Upsilon(T^{1}_{k+1})-R_{k})\Pi^{1*}_{k}-T^{1}_{k}\big] x_{k},x_{k}\big\rangle\\&\quad+2\big\langle\big[\Lambda(T^{1}_{k+1},\Pi^{2*}_{k})-L_{k}\big]^{\top}\Sigma^{1}_{k}+\Theta(T^{1}_{k+1},\Pi^{2*}_{k},\phi^{1}_{k+1})\\&\quad-q_{k},x_{k}\rangle+2\langle \Phi(T^{1}_{k+1},\Pi^{2*}_{k},\phi^{1}_{k+1})-\rho_k+\big(\Upsilon(T^{1}_{k+1})\\&\quad-R_{k}\big)\Sigma^{1}_{k},\Pi^{1*}_{k}x_{k}\big\rangle+\langle(\Upsilon(T^{1}_{k+1})-R_{k})\Sigma^{1}_{k},\Sigma^{1}_{k}\rangle\\&\quad+2\langle\Phi(T^{1}_{k+1},\Pi^{2*}_{k},\phi^{1}_{k+1})-\rho_k,\Sigma^{1}_{k}\rangle+\langle T^{1}_{k+1}\big(C_{k}\Sigma^{2*}_{k}	\end{split}
				\end{equation}
			\begin{equation}\label{eq:4.42}
		\begin{split}&\quad+b_{k}+\big(F_{k}\Sigma^{2*}_{k}+\sigma_{k}\big)\omega_{k}\big),C_{k}\Sigma^{2*}_{k}+b_{k}+\big(F_{k}\Sigma^{2*}_{k}+\sigma_{k}\big)\omega_{k}\rangle\\&\quad+2\langle\phi^{1}_{k+1},C_{k}\Sigma^{2*}_{k}+b_{k}+F_{k}\Sigma^{2*}_{k}\sigma_k\omega_{k}\rangle-2\langle\phi^{1}_{k},x_{k}\rangle\biggr\}.
			\end{split}
		\end{equation}
		Consider the cost functional \eqref{eq:40.2}, we have:
	\begin{equation}\label{eq:4.43}
		\begin{split}
			&\mathcal{J}_1(0, \xi, \Pi^{1*}\mathbf{x}+\Sigma^{1},\Pi^{2*}\mathbf{x}+\Sigma^{2*})-\frac{1}{2}\mathbb{E}\bigl\{\langle G_{N} x_{N}, x_{N}\rangle \\&+ 2\langle g_{N}, x_{N}\rangle\bigr\}
			\\=&\frac{1}{2} \mathbb{E}\bigg\{ \sum_{k=0}^{N-1} \bigg[\langle \big[Q_{k}+ L_{k}^{\top}\Pi^{1*}_{k}+ \Pi^{1*\top}_{k}L_{k}+\Pi^{1*\top}_{k}R_{k}\Pi^{1*}_{k}\big]\\
			&\times x_{k}, x_{k}\rangle  +2\langle  L_{k}^{\top}\Sigma^{1}_{k}+q_{k}, x_{k}\rangle+2\langle R_{k}\Sigma^{1}_{k}+\rho_{k}, \Pi^{1*}_{k}x_{k}\rangle \\
			& +\langle R_{k}\Sigma^{1}_{k}, \Sigma^{1}_{k}\rangle+2\langle \rho_{k}, \Sigma^{1}_{k}\rangle\bigg]\bigg\}.
		\end{split}
	\end{equation}
		Combining \eqref{eq:4.42}--\eqref{eq:4.43} with \eqref{eq:5.15}, \eqref{eq:5.18} under conditions \eqref{eq:5.17} and \eqref{eq:5.19}, we obtain:
		\begin{equation}\nonumber
			\begin{split}
				&\mathcal{J}_1(0, \xi, \Pi^{1*}\mathbf{x}+\Sigma^{1},\Pi^{2*}\mathbf{x}+\Sigma^{2*})-\frac{1}{2}\mathbb{E}\bigl\{\langle T^{1}_{0} x_{0}, x_{0}\rangle\\&+2\langle \phi^{1}_{0}, x_{0}\rangle \bigr\}
				\\=&\frac{1}{2} \mathbb{E}\sum_{k=0}^{N-1}\biggl\{2\big\langle\Lambda(T^{1}_{k+1},\Pi^{2*}_{k})^{\top}\Sigma^{1}_{k}+\Theta(T^{1}_{k+1},\Pi^{2*}_{k},\phi^{1}_{k+1}),x_{k}\rangle\\&+2\langle \Phi(T^{1}_{k+1},\Pi^{2*}_{k},\phi^{1}_{k+1})+\Upsilon(T^{1}_{k+1})\Sigma^{1}_{k},\Pi^{1*}_{k}x_{k}\big\rangle\\&+\langle\Upsilon(T^{1}_{k+1})\Sigma^{1}_{k},\Sigma^{1}_{k}\rangle+2\langle\Phi(T^{1}_{k+1},\Pi^{2*}_{k},\phi^{1}_{k+1}),\Sigma^{1}_{k}\rangle\\&+\langle T^{1}_{k+1}\big(C_{k}\Sigma^{2*}_{k}+b_{k}+\big(F_{k}\Sigma^{2*}_{k}+\sigma_{k}\big)\omega_{k}\big),C_{k}\Sigma^{2*}_{k}\\&+b_{k}+\big(F_{k}\Sigma^{2*}_{k}+\sigma_{k}\big)\omega_{k}\rangle+2\langle\phi^{1}_{k+1},C_{k}\Sigma^{2*}_{k}+b_{k}\\&+F_{k}\Sigma^{2*}_{k}\sigma_k\omega_{k}\rangle-2\langle\phi^{1}_{k},x_{k}\rangle\biggr\}
					\\=&\frac{1}{2} \mathbb{E}\sum_{k=0}^{N-1}\biggl\{2\big\langle\big(\Lambda(T^{1}_{k+1},\Pi^{2*}_{k})+\Upsilon(T^{1}_{k+1})\Pi^{1*}_{k}\big)^{\top}\Sigma^{1}_{k},x_{k}\rangle\\&+2\langle\Theta(T^{1}_{k+1},\Pi^{2*}_{k},\phi^{1}_{k+1}),x_{k}\rangle+2\langle \Phi(T^{1}_{k+1},\Pi^{2*}_{k},\phi^{1}_{k+1}),\\&\quad\Pi^{1*}_{k}x_{k}\big\rangle-2\langle\Lambda(T^{1}_{k+1},\Pi^{2*}_{k})^{\top}\Sigma^{1*}_k+\Theta(T^{1}_{k+1},\Pi^{2*}_{k},\phi^{1}_{k+1}),\\& x_{k}\rangle+\langle\Upsilon(T^{1}_{k+1})\Sigma^{1}_{k},\Sigma^{1}_{k}\rangle+2\langle\Phi(T^{1}_{k+1},\Pi^{2*}_{k},\phi^{1}_{k+1}),\Sigma^{1}_{k}\rangle\\&+\langle T^{1}_{k+1}\big(C_{k}\Sigma^{2*}_{k}+b_{k}+\big(F_{k}\Sigma^{2*}_{k}+\sigma_{k}\big)\omega_{k}\big),C_{k}\Sigma^{2*}_{k}+b_{k}\\&+\big(F_{k}\Sigma^{2*}_{k}+\sigma_{k}\big)\omega_{k}\rangle+2\langle\phi^{1}_{k+1},C_{k}\Sigma^{2*}_{k}+b_{k}+F_{k}\Sigma^{2*}_{k}\sigma_k\omega_{k}\rangle\biggr\}
				\\=&\frac{1}{2} \mathbb{E}\sum_{k=0}^{N-1}\biggl\{\langle\Upsilon(T^{1}_{k+1})\Sigma^{1}_{k},\Sigma^{1}_{k}\rangle-2\langle\Upsilon(T^{1}_{k+1})\Sigma^{1*}_{k},\Sigma^{1}_{k}\rangle\\&+\langle\Upsilon(T^{1}_{k+1})\Sigma^{1*}_{k},\Sigma^{1*}_{k}\rangle-\langle\Upsilon(T^{1}_{k+1})\Sigma^{1*}_{k},\Sigma^{1*}_{k}\rangle\\&+\langle T^{1}_{k+1}\big(C_{k}\Sigma^{2*}_{k}+b_{k}+\big(F_{k}\Sigma^{2*}_{k}+\sigma_{k}\big)\omega_{k}\big),C_{k}\Sigma^{2*}_{k}\\&+b_{k}+\big(F_{k}\Sigma^{2*}_{k}+\sigma_{k}\big)\omega_{k}\rangle+2\langle\phi^{1}_{k+1},C_{k}\Sigma^{2*}_{k}+b_{k}\\&+F_{k}\Sigma^{2*}_{k}\sigma_k\omega_{k}\rangle\biggr\}
				\\=&\frac{1}{2} \mathbb{E}\sum_{k=0}^{N-1}\biggl\{\langle \Upsilon(T^{1}_{k+1})\big(\Sigma^{1}_{k}-\Sigma^{1*}_{k}\big),\Sigma^{1}_{k}-\Sigma^{1*}_{k}\rangle-\langle\Upsilon(T^{1}_{k+1})\\&\times\Sigma^{1*}_{k},\Sigma^{1*}_{k}\rangle+\langle T^{1}_{k+1}\big(C_{k}\Sigma^{2*}_{k}+b_{k}+\big(F_{k}\Sigma^{2*}_{k}+\sigma_{k}\big)\omega_{k}\big),	\end{split}
			\end{equation}
		\begin{equation}\label{eq:4.44}
	\begin{split}&C_{k}\Sigma^{2*}_{k}+b_{k}+\big(F_{k}\Sigma^{2*}_{k}+\sigma_{k}\big)\omega_{k}\rangle+2\langle\phi^{1}_{k+1},C_{k}\Sigma^{2*}_{k}+b_{k}\\&+F_{k}\Sigma^{2*}_{k}\sigma_k\omega_{k}\rangle\biggr\}
						.
			\end{split}
		\end{equation}
	Consequently,
		\begin{equation}\nonumber
		\begin{split}
			&\mathcal{J}_1(0, \xi, \Pi^{1*}\mathbf{x}^{\ast}+\Sigma^{1}, \Pi^{2*}\mathbf{x}+\Sigma^{2\ast})\\&- \mathcal{J}_1(0, \xi, \Pi^{1*}\mathbf{x}+\Sigma^{1\ast}, \Pi^{2*}\mathbf{x}+\Sigma^{2\ast})\\=&\frac{1}{2} \mathbb{E}\sum_{k=0}^{N-1}\bigg\{\langle \Upsilon(T^{1}_{k+1})\big(\Sigma^{1}_{k}-\Sigma^{1*}_{k}\big),\Sigma^{1}_{k}-\Sigma^{1*}_{k}\rangle\bigg\}.
		\end{split}
	\end{equation}
			it follows that for any $\Sigma^{1}\in\mathcal{U}$,
			\begin{equation}\nonumber
				\begin{split}
					\mathcal{J}_1(0, \xi, \Pi^{1*}\mathbf{x}^{\ast}&+\Sigma^{1\ast}, \Pi^{2*}\mathbf{x}+\Sigma^{2\ast}) \\\leq&  \mathcal{J}_1(0, \xi, \Pi^{1*}\mathbf{x}+\Sigma^{1}, \Pi^{2*}\mathbf{x}+\Sigma^{2\ast}),
				\end{split}
			\end{equation}
			if and only if
			$$\Upsilon(T^{1}_{k+1}) \succeq 0, \quad k\in\mathcal{T}. $$
			Similarly, by Value function \eqref{eq:4.65}, for any $\Sigma^{2}\in\mathcal{V}$,
			\begin{equation}\nonumber
				\begin{split}
					\mathcal{J}_2(0, \xi, \Pi^{1*}\mathbf{x}^{\ast}&+\Sigma^{1\ast}, \Pi^{2*}\mathbf{x}+\Sigma^{2\ast}) \\\leq&  \mathcal{J}_2(0, \xi, \Pi^{1*}\mathbf{x}+\Sigma^{1\ast}, \Pi^{2*}\mathbf{x}+\Sigma^{2}),
				\end{split}
			\end{equation}
			if and only if
			$$\Upsilon(T^{2}_{k+1}) \succeq 0, \quad k\in\mathcal{T}.  $$
			This proves the sufficiency, as well as the necessity of \eqref{eq:5.16} and  \eqref{eq:5.21}.
		\end{proof}
		
		\section{Conclusion}\label{sec:5}
		
		This paper explores a class of LQ stochastic two-person non-zero-sum difference games with random coefficients, emphasizing the analysis of closed-loop Nash equilibria. The stochastic nature of the coefficients and weighting matrices in both the state equation and cost functionals leads to the emergence of fully coupled FBS$\Delta$Es, which are decoupled by two symmetric CCREs. Besides, we demonstrate that the existence of a closed-loop Nash equilibrium hinges on the solvability of these Riccati equations without strong regularity conditions. Furthermore, the closed-loop Value function is constructed naturally, providing a direct link between the equilibrium and the solution of the FBS$\Delta$Es.~Finally, the necessary and sufficient conditions for the closed-loop Nash equilibrium are characterized by a class of Lyapunov-type equations, obtained by the approach of DPP.~This work advances the theoretical framework for closed-loop Nash equilibria in stochastic LQ games, offering new perspectives for addressing such problems in future research.
		\bibliographystyle{plain}
		
		\section*{Appendix \\\\{A.~Simplified Notation for CCREs and Related Variables}}
		\addcontentsline{toc}{section}{Appendix A:~Simplified Notation for CCREs and Related Variables}
		\normalsize{To simplify the notation, the following symbols will be adopted.~For
$T^{1}, T^{2}   \in \mathbb{S}^{n}_{+},
\Pi^{1} \in \mathbb{R}^{m \times n},
\Pi^{2} \in \mathbb{R}^{l \times n},
\phi^{1}, \phi^{2} \in \mathbb{R}^{n},
\Sigma^{1} \in \mathbb{R}^{m},
\Sigma^{2} \in \mathbb{R}^{l}$, we define:
		\small{
			\begin{equation}\nonumber
				\begin{array}{l}
					\begin{split}
						\Upsilon&(T^{1}_{k+1}) \\\equiv& R_{k}+B^{\top}_{k}\mathbb{E}[T^{1}_{k+1} \mid \mathcal {F}_{k-1}]B_{k}+B^{\top}_{k}\mathbb{E}[T^{1}_{k+1}\omega_{k} \mid \mathcal {F}_{k-1}] E_{k}\\&+E^{\top}_{k}\mathbb{E}[T^{1}_{k+1}\omega_{k} \mid \mathcal {F}_{k-1}]B_{k}+E^{\top}_{k}\mathbb{E}[T^{1}_{k+1}\omega_{k}^{2} \mid \mathcal {F}_{k-1}]E_{k},
						\\
						\Lambda&\left(T^{1}_{k+1},\Pi^{2}_{k}\right) \\\equiv& L_{k}+B_{k}^{\top}\mathbb{E}[T^{1}_{k+1}\mid \mathcal {F}_{k-1}]\big(A_{k}+C_{k}\Pi^{2}_{k}\big)\\&+E_{k}^{\top}\mathbb{E}[T^{1}_{k+1}\omega_{k}\mid \mathcal {F}_{k-1}]\big(A_{k}+C_{k}\Pi^{2}_{k}\big)
						\\&+B_{k}^{\top}\mathbb{E}[T^{1}_{k+1}\omega_{k}\mid \mathcal {F}_{k-1}]\big(D_{k}+F_{k}\Pi^{2}_{k}\big)\\&+E_{k}^{\top}\mathbb{E}[T^{1}_{k+1}\omega_{k}^{2}\mid \mathcal {F}_{k-1}]\big(D_{k}+F_{k}\Pi^{2}_{k}\big),
						\\
						\Upsilon&(T^{2}_{k+1}) \\\equiv& S_{k}+C^{\top}_{k}\mathbb{E}[T^{2}_{k+1} \mid \mathcal {F}_{k-1}]C_{k}+C^{\top}_{k}\mathbb{E}[T^{2}_{k+1}\omega_{k} \mid \mathcal {F}_{k-1}] F_{k}\\&+F^{\top}_{k}\mathbb{E}[T^{2}_{k+1}\omega_{k} \mid \mathcal {F}_{k-1}]C_{k}+F^{\top}_{k}\mathbb{E}[T^{2}_{k+1}\omega_{k}^{2} \mid \mathcal {F}_{k-1}]F_{k},
				\\
						\Lambda&\left(T^{2}_{k+1},\Pi^{1}_{k}\right) \\\equiv&
						M_{k}+C_{k}^{\top}\mathbb{E}[T^{2}_{k+1}\mid \mathcal {F}_{k-1}]\big(A_{k}+B_{k}\Pi^{1}_{k}\big)\\&+C_{k}^{\top}\mathbb{E}[T^{2}_{k+1}\omega_{k}\mid \mathcal {F}_{k-1}]\big(D_{k}+E_{k}\Pi^{1}_{k}\big)\\&+F_{k}^{\top}\mathbb{E}[T^{2}_{k+1}\omega_{k}\mid \mathcal {F}_{k-1}]\big(A_{k}+B_{k}\Pi^{1}_{k}\big)\\&+F_{k}^{\top}\mathbb{E}[T^{2}_{k+1}\omega_{k}^{2}\mid \mathcal {F}_{k-1}]\big(D_{k}+E_{k}\Pi^{1}_{k}\big).
						\\
						\Delta&\left(T^{1}_{k+1},\Pi^{2}_{k}\right) \\\equiv& Q_{k}+\big(A_{k}+C_{k}\Pi^{2}_{k}\big)^{\top}\mathbb{E}[T^{1}_{k+1}\mid \mathcal {F}_{k-1}]\big(A_{k}+C_{k}\Pi^{2}_{k}\big)\\&+\big(A_{k}+C_{k}\Pi^{2}_{k}\big)^{\top}\mathbb{E}[T^{1}_{k+1}\omega_{k}\mid \mathcal {F}_{k-1}]\big(D_{k}+F_{k}\Pi^{2}_{k}\big)\\&+\big(D_{k}+F_{k}\Pi^{2}_{k}\big)^{\top}\mathbb{E}[T^{1}_{k+1}\omega_{k}\mid \mathcal {F}_{k-1}]\big(A_{k}+C_{k}\Pi^{2}_{k}\big)
						\\&+\big(D_{k}+F_{k}\Pi^{2}_{k}\big)^{\top}\mathbb{E}[T^{1}_{k+1}\omega_{k}^{2}\mid \mathcal {F}_{k-1}]\big(D_{k}+F_{k}\Pi^{2}_{k}\big),
				\\
						\Delta&\left(T^{2}_{k+1},\Pi^{1}_{k}\right)\\\equiv&P_{k}+\big(A_{k}+B_{k}\Pi^{1}_{k}\big)^{\top}\mathbb{E}[T^{2}_{k+1}\mid \mathcal {F}_{k-1}]\big(A_{k}+B_{k}\Pi^{1}_{k}\big)\\&+\big(A_{k}+B_{k}\Pi^{1}_{k}\big)^{\top}\mathbb{E}[T^{2}_{k+1}\omega_{k}\mid \mathcal {F}_{k-1}]\big(D_{k}+E_{k}\Pi^{1}_{k}\big)\\&+\big(D_{k}+E_{k}\Pi^{1}_{k}\big)^{\top}\mathbb{E}[T^{2}_{k+1}\omega_{k}\mid \mathcal {F}_{k-1}]\big(A_{k}+B_{k}\Pi^{1}_{k}\big)
						\\&+\big(D_{k}+E_{k}\Pi^{1}_{k}\big)^{\top}\mathbb{E}[T^{2}_{k+1}\omega_{k}^{2}\mid \mathcal {F}_{k-1}]\big(D_{k}+E_{k}\Pi^{1}_{k}\big)
						\\
						\Phi&\left(T^{1}_{k+1},\Pi^{2}_{k},\phi^{1}_{k+1}\right)\\\equiv&\rho_{k}+B^{\top}_{k}\bigg[\mathbb{E}[T^{1}_{k+1}\mid \mathcal {F}_{k-1}]\big(C_{k}\Sigma^{2}_{k}+b_{k}\big)+\mathbb{E}[T^{1}_{k+1}\omega_{k}\mid \mathcal {F}_{k-1}]\\&\big(F_{k}\Sigma^{2}_{k}+\sigma_{k}\big)+\mathbb{E}[\phi^{1}_{k+1}\mid \mathcal {F}_{k-1}]\bigg]+E^{\top}_{k}\bigg[\mathbb{E}[T^{1}_{k+1}\omega_{k}\mid \mathcal {F}_{k-1}]\\&\big(C_{k}\Sigma^{2}_{k}+b_{k}\big)+\mathbb{E}[T^{1}_{k+1}\omega_{k}^{2}\mid \mathcal {F}_{k-1}]\big(F_{k}\Sigma^{2}_{k}+\sigma_{k}\big)\\&+\mathbb{E}[\phi^{1}_{k+1}\omega_{k}\mid \mathcal {F}_{k-1}]\bigg]
						\\
						\Phi&\left(T^{2}_{k+1},\Pi^{1}_{k},\phi^{2}_{k+1}\right)\\\equiv&\theta_{k}+C^{\top}_{k}\bigg[\mathbb{E}[T^{2}_{k+1}\mid \mathcal {F}_{k-1}]\big(B_{k}\Sigma^{1}_{k}+b_{k}\big)+\mathbb{E}[T^{2}_{k+1}\omega_{k}\mid \mathcal {F}_{k-1}]\\&\big(E_{k}\Sigma^{1}_{k}+\sigma_{k}\big)+\mathbb{E}[\phi^{2}_{k+1}\mid \mathcal {F}_{k-1}]\bigg]+F^{\top}_{k}\bigg[\mathbb{E}[T^{2}_{k+1}\omega_{k}\mid \mathcal {F}_{k-1}]\\&\big(B_{k}\Sigma^{1}_{k}+b_{k}\big)+\mathbb{E}[T^{2}_{k+1}\omega_{k}^{2}\mid \mathcal {F}_{k-1}]\big(E_{k}\Sigma^{1}_{k}+\sigma_{k}\big)\\&+\mathbb{E}[\phi^{2}_{k+1}\omega_{k}\mid \mathcal {F}_{k-1}]\bigg],
						\\\Theta&\left(T^{1}_{k+1},\Pi^{2}_{k},\phi^{1}_{k+1}\right) \\\equiv& q_{k}+\big(A_{k}+C_{k}\Pi^{2}_{k}\big)^{\top}\bigg[\mathbb{E}[T^{1}_{k+1}\mid \mathcal {F}_{k-1}]\big(C_{k}\Sigma^{2}_{k}+b_{k}\big)\\&
						+\mathbb{E}[T^{1}_{k+1}\omega_{k}\mid \mathcal {F}_{k-1}]\big(F_{k}\Sigma^{2}_{k}+\sigma_{k}\big)+\mathbb{E}[\phi^{1}_{k+1}\mid \mathcal {F}_{k-1}]\bigg]
					\end{split}
			\end{array}
		\end{equation}}
	\footnotesize{\begin{equation}\label{eq:2.30}
\begin{array}{l}
\begin{split}
						&+\big(D_{k}+F_{k}\Pi^{2}_{k}\big)^{\top}\bigg[\mathbb{E}[T^{1}_{k+1}\omega_{k}\mid \mathcal {F}_{k-1}]\big(C_{k}\Sigma^{2}_{k}+b_{k}\big)\\&+\mathbb{E}[T^{1}_{k+1}\omega_{k}^{2}\mid \mathcal {F}_{k-1}]\big(F_{k}\Sigma^{2}_{k}+\sigma_{k}\big)
						+\mathbb{E}[\phi^{1}_{k+1}\omega_{k}\mid \mathcal {F}_{k-1}]\bigg],
						\\
						\Theta&\left(T^{2}_{k+1},\Pi^{1}_{k},\phi^{2}_{k+1}\right) \\\equiv& p_{k}+\big(A_{k}+B_{k}\Pi^{1}_{k}\big)^{\top}	\bigg[\mathbb{E}[T^{2}_{k+1}\mid \mathcal {F}_{k-1}]\big(B_{k}\Sigma^{1}_{k}+b_{k}\big)\\&
						+\mathbb{E}[T^{2}_{k+1}\omega_{k}\mid \mathcal {F}_{k-1}]\big(E_{k}\Sigma^{1}_{k}+\sigma_{k}\big)+\mathbb{E}[\phi^{2}_{k+1}\mid \mathcal {F}_{k-1}]\bigg]
						\\&+\big(D_{k}+E_{k}\Pi^{1}_{k}\big)^{\top}\bigg[\mathbb{E}[T^{2}_{k+1}\omega_{k}\mid \mathcal {F}_{k-1}]\big(B_{k}\Sigma^{1}_{k}+b_{k}\big)\\&+\mathbb{E}[T^{2}_{k+1}\omega_{k}^{2}\mid \mathcal {F}_{k-1}]\big(E_{k}\Sigma^{1}_{k}+\sigma_{k}\big)
						+\mathbb{E}[\phi^{2}_{k+1}\omega_{k}\mid \mathcal {F}_{k-1}]\bigg],
					\end{split}
				\end{array}
		\end{equation}}
		
	\end{document}